   \def\obs#1{{\bf (*** #1 ***)} }

 \def\obs#1{}     

\documentclass[twoside,letterpaper,draft,11pt]{amsart}

\usepackage{amsmath}               
\usepackage{amsthm}                
\usepackage{latexsym}

\usepackage{amssymb}
\usepackage{xspace}
\usepackage{amscd}

\usepackage{mathrsfs}

\usepackage[all]{xy}

\theoremstyle{definition}

\theoremstyle{remark}

\newcommand{\bbC}{{\mathbb C}}

\newcommand{\D}{{\mathcal D}}

 \newcommand{\af}{\alpha}

\newcommand{\Q}{{\mathbb Q}}

\newcommand{\A}{{\mathcal A}}

\newcommand{\B}{{\mathcal B}}

\newcommand{\V}{{\mathcal V}}
\newcommand{\U}{{\mathcal U}}

\newcommand{\C}{{\mathbb C}}           
\newcommand{\Z}{{\mathbb Z}}           
   
\newcommand{\N}{{\mathbb N}}             
\newcommand{\F}{{\mathbb F}}           
\newcommand{\kpar}{{\kappa}_{\rm par}}     
\newcommand{\Cpar}{\C_{\rm par}}

\newcommand{\m}{{}^{-1}}

\newcommand{\cG}[1]%
{%
	\mathcal G{(#1)}%
}

\newcommand{\tq}{\;\mid\;}

\newcommand{\lpar}{\lambda_{\rm p}} 

\begin{document}

\title[Recent  developments around partial actions]{Recent  developments around partial actions}

\author{M. Dokuchaev}
\address{Instituto de
Matem\'atica e Estat\'\i stica\\
Universidade de S\~ao Paulo\\
05508-090 S\~ao Paulo, SP, Brasil}
\email{dokucha@gmail.com}
\thanks{This work was partially supported by CNPq and Fapesp of Brazil}
\keywords{Partial action, partial representation, crossed product.}

\dedicatory{Dedicated to Antonio Paques on the occasion of his 70th birthday}


\subjclass[2000]{Primary 16W22,  16S35, 20C99, 46L55; Secondary 08A02, 13B05,  16S10,   16S36, 16T05, 16T15, 16W50, 20C15, 20C25, 20F10, 20L05, 20M18,  20M25, 20M30, 20M50, 22A22, 46L05, 54H15.}

\date{}

\begin{abstract} We give an overview of publications on partial actions and related concepts, paying main attention to some recent developments.
\end{abstract}

\maketitle


Since the publication of our previous survey \cite{D2} the ac\-ti\-vi\-ty around  partial actions enjoyed an 
increasing intensity, resulting in  re\-mar\-kab\-le 
applications and theoretic development, as well as the  publication  of the book \cite{E6} by Ruy~Exel. The latter   contains  algebraic and  $C^*$-algebraic background, a comprehensive treatment of 
graded $C^*$-algebras via Fell bundles and partial $C^*$-crossed products, as well as  prominent applications to 
Wiener-Hopf $C^*$-algebras associated to  quasi-lattice ordered groups and graph $C^*$-algebras.    The working 
team gained new participants, in particular, a number of young researchers have been brought to the area by their supervisors.  Dedicating this 
survey to  Antonio~Paques on the occasion of his  $70$th birthday,  I take the opportunity  to register his remarkable 
leadership in introducing junior scientists  into the subject.\\

For a reader not familiar with  our topic, we recall that partial actions appear naturally restricting usual (global)  actions as follows. Let $\beta $  be a global  action of a 
group $G$ on a set $Y,$ i.e. we have a homomorphism $\beta : G \ni g \mapsto \beta_g \in S(Y),$ from $G$ into the symmetric group 
$S(Y)$ of all bijections (permutations) $Y\to Y.$ Now, let $X\subseteq Y$ be a subset. For each $g\in G$ denote 
\begin{equation}\label{restriction}
X_g = \beta_g (X) \cap X,
\end{equation}
  and restricting $ \beta_g $ to $X_{g\m}$ we obtain a bijection   $\alpha _g : X_{g\m} \to X_g$ between subsets of $X$ (partial bijections of $X$). Then the collection  $\alpha = \{   \alpha _g \, : \,  g\in G  \}$ gives us  a partial action of $G$ on $X.$ Notice that obviously\\ 

\noindent (i) $\alpha _1 = {\rm id}_X,$\\

\noindent and for all $x \in X, g,h \in G:$\\

\noindent (ii) $ \exists \; \alpha_h(x),  \; \exists \; \alpha_g (\alpha_h(x)) \; \Longrightarrow   \;   \exists  \; \alpha _{gh}(x) \;
\; \text{and} \;\;    \alpha _g (\alpha_h(x)) = \alpha _{gh}(x).$\\

\noindent Thus one defines a {\it partial action} $\alpha $ of a group $G$ on a set $X$ as a family of partial bijections 
$\alpha _g : X_{g\m} \to X_g, $ $  (g\in G) $ of $X,$ such that (i) and (ii) are sa\-tis\-fi\-ed. The subsets $X_g$ are called the 
domains of $\alpha $, and one refers to the triple $(X, G, \alpha $) (or, more precisely, to the quadruple 
$(X, G, \{ X_g \}_{g\in G},  \{ \alpha _g \}_{g\in G})$) as a {\it partial dynamical system}.\\ 

If $X$ has some structure, then one imposes appropriate restrictions on the $X_g$ and $\alpha_g.$ In particular, if we 
define a partial action on a  ring or an algebra (or, more generally, a multiplicative semigroup), then we assume that each domain  is a two-sided ideal  and each 
$\alpha _g$ is an isomorphism 
of rings (algebras, or semigroups).\footnote{Notice that in \cite{CaenJan} the authors defined a partial group action  on an algebra assuming that each $X_g$ is a right ideal generated by an idempotent.}   If our ring has an involution then the domains should be closed under the involution and each $\alpha _g$ is supposed to respect the involution.   If $X$ is a topological space and $G$ is a discrete group, then we assume that  each $X_g$ is an open subset and  each $\alpha _g$ is a homeomorphism.  In this case we say that  $(X, G, \alpha $) is a {\it topological partial dynamical system}. The latter term   englobes also  a partial action of a  topological group $G$ on a topological space, for the definition of which   more requirements are needed 
(see \cite[Definition 1.1]{AbadieTwo}). When definining a partial action of a discrete group $G$ on a $C^*$-algebra one 
assumes that each domain is a norm closed ideal and each $\alpha _g$ is a $*$-isomorphism of  algebras. 
In this case we have a {\it $C^*$-algebraic  partial dynamical system}. The more general 
definition of a (continuous) partial action of a locally compact group  on a   $C^*$-algebra was given  by R.~Exel in \cite{E0}.\\

A partial  action of a group $G$ on an algebra gives rise to two more concepts: that of a partial action of $G$ on a (multiplicative) semigroup (or monoid) and the notion of a partial action of $G$ on a vector space (or, more ge\-ne\-ral\-ly, a $\kappa$-module over a commutative ring $\kappa $). In the first case  each domain $X_g$ is assumed to be a two-sided ideal in the semigroup  and every $\alpha _g$ is a semigroup isomorphism. In the second case  the $X_g$ are subspaces ($\kappa$-submodules) and the $\alpha _g$ are $\kappa $-linear isomorphisms.\\

As we saw above, any global group action on a set $Y$ can be restricted to any subset $X,$ resulting in a partial action whose domains are defined by (\ref{restriction}). With some structure on $X$  the above procedure works with an appropriate  assumption on $Y.$ Thus any  group action on a topological space can be restricted to an open subset, and any group action on a ring restricts this way to an arbitrary (two-sided) ideal. This suggests  an important question: is a given partial action can be seen (after some identifications) as a restriction of a global one? In the case of the partial actions on sets the 
answer is always positive \cite{AbadieTwo} (see also \cite[Theorem 3.5]{E6}), however,  as we shall see below, the situation is more complex with the presence of  some structure.\\

   The definition of a partial action of a group $G$ on a set $X$ can be reformulated in several ways. Denote by 
${\mathcal I}(X)$  the set of all partial bijections of $X,$ including the vacuous  bijection $\emptyset \to \emptyset.$ One defines the composition of two  partial bijections $\varphi$ and $\psi$ of $X$ on the largest possible subset: 

$$ \phi \circ \psi : \;\; {\psi}\m ( {\rm ran}\, \psi  \cap  {\rm dom}\, \varphi )  \to \varphi (  {\rm ran}\, \psi  \cap  {\rm dom}\, \varphi ).$$ 

\noindent The above operation endows ${\mathcal I}(X)$ with a structure of an inverse monoid with zero  $\emptyset \to \emptyset ,$ called the symmetric inverse semigroup of $X.$ Then a partial action $\alpha $ of $G$ on $X$ gives a map $\alpha : G \ni g \mapsto \alpha_g \in   {\mathcal I}(X),$ and one may wonder 
which maps  $G \to  {\mathcal I}(X)$ correspond to partial actions on $X.$\\

 The first answer was given by R.~Exel in 
\cite[Proposition 4.1]{E1} relating it  to the important concept of a partial representation:  a map 
$\alpha : G \to    {\mathcal I}(X)$ gives a partial action if and only if for all $g,h\in G$ we have

\hspace*{5mm}

\noindent  (a) $\alpha_1 =  {\rm id}_X,$

\hspace*{1mm}

\noindent  (b) $\alpha _{g\m}   \alpha _g \alpha _h = \alpha _{g\m}   \alpha _{gh}.$

\hspace*{1mm}

\noindent  In this case $\alpha$ also satisfies 

\hspace*{1mm}

\noindent    (c) $\alpha _{g}   \alpha _h \alpha _{h\m} = \alpha _{g h}   \alpha _{h\m}.$\\

The passage between (b) and (c) can be easily performed by applying the inverses of the partial bijections in  
$ {\mathcal I}(X).$ This  result leads to the following relevant notion. A map $\alpha :  G \to S$ from $G$ to a monoid  $S$ is said to be a (unital) {\it partial homomorphism} if $\alpha $ sends $1_G$ to $1_S$ and   satisfies  (b) and (c) above.  If $S$ is the multiplicative monoid  of a unital  algebra, then we say that $\alpha $ is a {\it partial representation}. This is the most basic relation between partial actions and partial representations, but there are others resulting in  a fruitful interaction, both theoretical and practical.\\

The inverse semigroup $ {\mathcal I}(X)$ possesses a natural partial order given by restriction of partial bijections. More generally,    
every inverse semigroup $S$ has a natural partial order defined by setting $a\leq b,$ $(a,b \in S)$  if and only if there exists 
an idempotent $e\in S$ with $a =be$ (see, for example, \cite{L}). It is immediately seen that (ii) above means that the composition $\alpha _g \circ \alpha _h$ is a restriction of $\alpha _{gh},$
so that one may write $\alpha _g \circ \alpha _h \leq \alpha _{gh},$ which is clearly a defining property of partial actions.
One may easily  obtain another one, i.e. $\alpha _{g\m} = ({\alpha} _g)\m.$ 
Taking this into account,   J. Kellendonk and M. V. Lawson \cite{KL1}, \cite{KL} gave another way to characterize   the maps   $G \to    {\mathcal I}(X)$ which determine  partial actions. To formalize this, recall that a function $\alpha : S \to T$ between inverse monoids  is called a {\it unital  premorphism} \cite{KL} (also referred to as   a {\it dual  prehomomorphism} \cite{LawMarSte}) if  for all $g,h \in S$ the following properties hold:

\hspace*{3mm}

\noindent  (1) $\alpha (1_S)  =  1_T,$

\hspace*{1mm}

\noindent  (2) $\alpha (g\m)  = \alpha (g)\m  ,$

\hspace*{1mm}

\noindent    (3) $\alpha (g)   \alpha (h) \leq \alpha (g h)  .$\\

\noindent Then a map $\alpha : G \to     {\mathcal I}(X)$ determines a partial action precisely when $\alpha $ is a unital premorphism.\\ 

Here we warn the reader with respect to the use of terminology. By a partial action of a group $G$ on a set $X$ the authors of 
\cite{KL} mean a family of partial bijections     $\alpha = \{   \alpha _g : X_{g\m} \to X_g \, : \,  g\in G  \}$ such that the corresponding map $\alpha : G \to     {\mathcal I}(X)$   is a non-necessarily unital premorphism, i.e. $\alpha $ satisfies  conditions (2) and (3).  If $\alpha $ also satisfies (1), then $\alpha $ is called  unital in \cite{KL}. Nevertheless, we prefer to assume  that (1) is always satisfied, and use the term unital in a different sense.  Namely, by a {\it unital partial action} we mean a partial action of group $G$ on a ring $\A $ 
such that each domain is a unital ring, i.e. generated by an idempotent which is central in $\A .$ These are exactly those partial actions of $G$ on a unital ring $\A$ which are glo\-ba\-li\-za\-b\-le (see \cite{DE}). Moreover, if each domain is an $s$-unital ring, then $\alpha $ is called an  {\it $s$-unital  partial action} on $\A $ (see \cite{AbDoExSi} and \cite{BPi}). Such partial actions appeared in the study of the globalization problem for partial group actions on $s$-unital rings in \cite{DdRS}.\\

 One should also notice that the term ``partial action'' of a group was also used in the literature in a close but different sense, and we refer the reader to the  survey \cite{D2} for details.\\

As it is was mentioned in many occasions, the formal concept of a partial action (in the sense we are using it)  appeared  in the  
theory $C^*$-algebras (see \cite{E-1}, \cite{Mc}, \cite{E0}, \cite{E1}), permitting one to endow relevant classes of   $C^*$-algebras with a general structure of a partial crossed product  \cite{E-2},  \cite{E-3},   \cite {E2}, \cite{ELQ}, \cite{QR} (see also \cite{E6}), 
and promptly  stimulating further use and discussions  in the area  \cite{Abadie}, \cite{AbadieTwo},    
\cite{AEE},  \cite {E2}, \cite{E2.5},  \cite{EL}, \cite{EL1},  \cite{QR}, \cite{Sieben}, \cite{Sieben98}, \cite{Sieben2001}. Subsequent    $C^*$-algebraic and topological developments on partial actions were made in \cite{AbadieThree},   \cite{EL2},  \cite{Megre2},  \cite{AbadieFour}, \cite{DonsHopen}, \cite{Lebedev}, \cite{Boava},  \cite{CrispLaca}, \cite{ChLi1}, 
 \cite{AbadiePerez} and \cite{EV}.\\

 More recent $C^*$ and topological advances  include the groupoid approach to the enveloping  $C^*$-algebras associated to  partial actions of  countable discrete groups on  (locally) compact spaces in \cite{EGG}, the use of inverse semigroup expansions to treat 
$C^*$-crossed products by twisted partial actions via twisted global actions of the expansion in 
\cite{BuE},
 the full (respectively, reduced) partial $C^*$-crossed product descriptions of  full (respectively, reduced)
$C^*$-algebras of countable  $E$-unitary or strongly $0$-$E$-unitary inverse semigroups as well as of  tight groupoids of countable  strongly $0$-$E$-unitary inverse semigroups in \cite{MilanSt},
the  study of the continuous orbit equivalence for partial dynamical systems and of the partial transformation groupoids with applications to graph $C^*$-algebras and semigroup $C^*$-algebras  in \cite{XinLi}, the partial group action approach to produce  a Bratteli-Vershik model linked to a minimal homeomorphism between open subsets with finite disjoint complements of the Cantor set in \cite{GiordanoGonStarl}, new developments on the globalization problem for partial actions on $C^*$-algebras and Hilbert bimodules in \cite{Ferraro},
the employment of the partial crossed pro\-duct theory to the investigation of  the Cuntz-Li $C^*$-algebras  related  to an integral domain in \cite{BoE1} with a further development in \cite{2017Vieira} for $C^*$-algebras associated with an injective endomorphism of a group with finite cokernel. Moreover, partial crossed products turned out to be useful   to deal with   $C^*$-algebras ari\-sing from self-similar graph actions in \cite{EStar1}, with $C^*$-algebras associated to any  stationary ordered Bratteli diagram in \cite{GR0}, as well as with  ultragraph $C^*$-algebras and related infinite alphabet shifts in \cite{GR2}, \cite{GR3}, \cite{GR4} and \cite{GilesGon}. In addition, in    \cite{KrakenQuastTimm},   partial coactions of $C^*$-bialgebras, in par\-ti\-cu\-lar of $C^*$-quantum groups,    on $C^*$-algebras were defined and studied, and a globalization result was obtained.\\

   Partial dynamical systems  were  also applied in \cite{BAbadieFAbadie} to the study of the ideal structure of  full and reduced cross-sectional $C^*$-algebras of Fell bundles,  in \cite{AraEKa}   to the investigation  
 of    $C^*$-algebras of  dynamical systems of type $(m,n)$, in \cite{CarlLars2016} to relative graph $C^*$-algebras,       in \cite{DE2} to $C^*$-algebras associated to arbitrary subshifts,   in 
\cite{GiordanoSierak} to the  study of the ideals and the pure infiniteness of partial   $C^*$-crossed pro\-ducts, and in 
 \cite{RenWil} to     the investigation of groupoids arising from partial semigroup actions and topological higher rank graphs.  Furthermore, partial representations and cohomology based on partial actions turned out to be useful in dealing with  the  separation and intersection properties of  ideals in global reduced  $C^*$-crossed products  \cite{KennedySchafhauser}.
In addition,  continuous partial actions of Polish groups on Polish spaces, and more generally,  of separable metrizable groups on Hausdorff (in the majority of facts metrizable) spaces were studied in \cite{GoPiUz1}, \cite{GoPiUz2},  \cite{PiUz1} and  \cite{PiUz2}. In particular, Effros' theorem \cite{Effros} on orbits in Polish spaces under continuous actions of Polish groups is extended in  \cite{GoPiUz1} to the context of partial actions.\\

A remarkable recent application was achieved to paradoxical decompositions in  \cite{AraE1}, where for a finite bipartite separated graph $(E,C)$ the tame graph $C^*$-algebra ${\mathcal O}(E,C)$ was introduced and proved to be isomorphic to a partial crossed product of a commutative $C^*$-algebra by a finitely generated free group. This permitted the authors of    \cite{AraE1} to give a negative answer to  an open  problem   on paradoxical decompositions in a topological setting, posed in \cite{Kerr}, \cite{KerrNowak} and \cite{RordamSier}.\\

The first algebraic results on the subject appeared in \cite{E1}, \cite{DEP},  \cite{KL1}, \cite{St1}, \cite{ChLi}, \cite{St2}, \cite{KL},  \cite{Choi1} and \cite{DE}, including  the first algebraic application for  tiling semigroups in \cite{KL1} (with further use in \cite{Zhu}),  for $E$-unitary inverse semigroups in \cite{KL}, to inverse semigroups and $F$-inverse monoids in   \cite{St1} and to inverse monoids of M\"obius type in \cite{ChLi}.  Independently from Exel's definition of a partial action, T.  Coulbois \cite{Coul} used a more restrictive notion, called a pre-action, to deal with the  Ribes-Zalesski  property  $(RZ_n)$ of groups from model 
theoretic point of  view (see also \cite{Coul2}). Furthermore,   J.-C. Birget \cite{Birget} applied  the partial action definition of 
 Thompson's groups to study algorithmic problems for them.
Since then the algebraic approach is being developed in diverse directions in various levels of generality,
 in\-clu\-ding partial actions of 
Hopf (or, more generally,  weak Hopf) algebras  \cite{AAB}, \cite{AB1}, \cite{AB2}, \cite{AB3}, \cite{ABDP}, \cite{ABDP2}, 
\cite{ABV}, 
\cite{BatiVerc1},
\cite{CaenDGr}, \cite{CaenJan}, \cite{CasPaqQuaSant},  \cite{CasPaqQuaSant2},
\cite{CasQua},  \cite{CasQua2}, \cite{CavalStAna1},   \cite{CavalStAna2}, 
\cite{ChenWang2013}, \cite{ChenWang2014}, \cite{ChenWangKang2017},
\cite{FerGonRod}, \cite{Paques2016},
 \cite{Zhang2016},
semigroups   
\cite{BuE}, \cite{CornGould},  \cite{DNZholt}, \cite{GouHol1}, 
\cite{GouHol2}, \cite{Hol1},    \cite{Khry1}  \cite{Kud},  \cite{LawMarSte}, \cite{Limin},  \cite{Megre2}, \cite{MikhNovZholt}, \cite{NovPerZholt}, \cite{PeZholt}, 
inductive
constellations \cite{GouHol2},   
groupoids  \cite{B},  \cite{BFP}, \cite{BP}, \cite{BPi},  \cite{BPi2}, \cite{Gil}, and, more generally, categories \cite{Nystedt}.  
In particular,  further algebraic applications have been found to graded algebras \cite{DE}, \cite{DES1}, to Hecke algebras \cite{E3}, to Leavitt path algebras \cite{GonOinRoy}, \cite{GR}, \cite{GonYon}, \cite{NysOinPin}, to inverse semigroups  \cite{Khry1}, \cite{Limin}, to restriction semigroups \cite{CornGould}, \cite{Kud},  to automata and machines \cite{DNZholt}, \cite{MikhNovZholt},  \cite{NovPerZholt},  \cite{PeZholt}, and  to Steinberg algebras \cite{BeuGon2}, \cite{BeuGonOinRoy}.\\

At this point we observe that the term ``partial action'' of a monoid (or  semigroup) is used in two slightly different senses. First M.  G. Megrelishvili and  L. Schr\"oder in   \cite{Megre2} gave the following definition of a partial monoid action, extending   the notion of a partial group action:
given a monoid  $M$  with unit $e$ and a set  $X,$ a  {\it (left) partial action} of $M$ on $X$ is   a partially defined map $M\times X\to X,$ 
$(m,x) \mapsto m\cdot x,$ such that, for all $x \in X, u, v \in M:$\\

\noindent ($i$) $e \cdot x  = x,$ for all $x\in X,$\\

\noindent ($ii$) $ \exists \; v\cdot x,   \; \exists \; u \cdot (v\cdot x ) \; \Longrightarrow   \;   \exists  \; (uv) \cdot x \;
\; \text{and} \;\;    u \cdot (v\cdot x ) =  (uv) \cdot x,$\\

\noindent ($iii$)  $ \exists \; v\cdot x,   \; \exists \; (u v) \cdot x \; \Longrightarrow   \;   \exists  \;     u \cdot (v\cdot x) \;
\; \text{and} \;\;    u \cdot (v\cdot x) =  (uv) \cdot x.$\\ 

\noindent A right partial action is defined symmetrically. The above was called a {\it strong partial action} of a monoid in the article 
\cite{Hol1} by C.\ Hollings, whereas in his definition of a partial monoid action (sometimes also referred to as a {\it weak partial action})  only ($i$) and ($ii$) are assumed. Of course, if in the definition by C.\ Hollings (i.e. in ($i$) and ($ii$)), we replace $M$ by a group $G,$ we  obtain the concept of a partial group action in the form spelled out by  J. Kellendonk and 
M. V. Lawson in  \cite{KL}. Ignoring above ($i$) we obtain the notion of a strong partial action of a semigroup, and omitting both ($i$) ad ($iii$) we come to that of a (weak) partial semigroup action. The partial monoid actions considered \cite{CornGould}, \cite{Kud} are weak, whereas in \cite{GouHol1} both weak and strong partial monoid (or semigroup) actions are discussed. 
Strong partial monoid actions were used in  \cite{DNZholt} under the name ``preactions'' to define the concept of a {\it preautomaton} (see also \cite{MikhNovZholt},  \cite{NovPerZholt},  \cite{PeZholt}). Notice that the term  ``partial automaton'' is already in use in a larger sense. Observe that there is a third kind of ``partial actions'' of monoids (or semigroups) called in  \cite[Definition 3.2]{GouHol1}  {\it incomplete actions} (see also \cite[p. 297]{Hol1} and \cite[\S 3]{Renshaw}). The definition of an incomplete action is obtained by replacing ($ii$) and ($iii$) above by the following stronger requirement:\\

\noindent ($ii'$) $ \exists \; v\cdot x,   \; \exists \; u \cdot (v\cdot x ) \; \Longleftrightarrow   \;   \exists  \; (uv) \cdot x, \;
\; \text{in which case } \;\;    u \cdot (v\cdot x ) =  (uv) \cdot x,$\\  

\noindent $x \in X, u, v \in M.$ Each incomplete action is a strong partial action, but there are strong partial actions which are not incomplete (see  Examples 2 and 3 in  \cite{GouHol1}).  Note, in addition, that the class of weak partial actions also properly contains the class of strong partial actions \cite[Example 2.11]{Hol1}.\\

In the definition of a partial action of an inverse semigroup the language of premorphisms is usually adopted. More precisely, in      \cite{LawMarSte} and
\cite{BuE}  a partial action of an inverse semigroup $G$ on a set $X$ is defined as an order-preserving premorphism 
$G \to {\mathcal I}(X), $ whereas in  \cite{Khry1} an arbitrary premorphism $G \to {\mathcal I}(X) $ is considered as a partial action.  In \cite[Definition 2.11]{BuE} the notion of a partial homomorphism  of an inverse  semigroup $G$ into a  semigroup  is  defined and it is proved that a map  between inverse semigroups  is a partial homomorphism if and only if 
it is  an order-preserving premorphism \cite[Proposition 3.1]{BuE}. Thus, similarly to the group case, partial actions of inverse semigroups on sets from \cite{LawMarSte} and
\cite{BuE} can be characterized as  partial homomorphisms  of the form 
$G \to {\mathcal I}(X).$ Notice that in \cite{Limin} partial actions of primitive inverse semigroups on sets were considered and applied to $E^*$-unitary categorical inverse semigroups.\\

More generally, weak and strong right partial actions of a left restriction semigroup $S$ on a set $X$ are defined by  V. Gould and C. Hollings   in  \cite[Definitions 3.3, 3.4]{GouHol1}, by adding to the right-handed forms   of ($ii$) and ($iii$)  an additional axiom which takes into  account the presence of the unary ope\-ra\-tion $u \mapsto u^+:$\\

\noindent ($iv$)   $  \exists \; x \cdot u   \;  \Longrightarrow   \;   \exists  \;    x \cdot  u^+  \;
\; \text{and} \;\;   x \cdot  u ^+ =   x,$\\ 

\noindent $u\in S, x \in X.$   Left restriction semigroups may be viewed as an  axi\-o\-ma\-ti\-za\-tion of the   partial transformation monoid 
${\mathcal P}{\mathcal I}(X)$ (see \cite{GouHol1}, \cite{Gou}, \cite{Hol2}).  The latter   is defined as the  semigroup of all partial maps of $X,$ i.e. non-necessarily bijective  functions of the form $X\supseteq Y \to Z \subseteq X.$  Notice that each left restriction semigroup can be seen as a subsemigroup of the partial transformation monoid 
${\mathcal P}{\mathcal I}(X)$ of  some set $X,$  which is closed under the unary operation of ``taking domains'', i.e. 
$^+: \varphi \mapsto {\rm id}_{{\rm dom} \; \varphi }$ (see \cite[Theorem 2.2]{GouHol1}, \cite[Corollary 6.3]{Gou}). 
Since every inverse semigroup $G$ is isomorphic to a $^+$-closed subsemigroup of ${\mathcal I}(X) \subset {\mathcal P}{\mathcal I}(X),$ it follows that  $G$ is a left restriction semigroup. Another example of a left restriction semigroup is given by the Szendrei expansion of a monoid \cite[Proposition 3.3]{Hol2}.
It was shown in \cite[p. 367]{GouHol1} that  partial actions of an inverse semigroup $G,$ in the  sense of \cite{LawMarSte} and
\cite{BuE}, are exactly the strong partial actions of $G$ considered as a left restriction semigroup.\\

One of the relevant problems in the theory is that of the globalization: given a group 
(semigroup, groupoid etc.) $G$ acting partially on an object $X,$ construct an embedding of $X$ into a larger object $Y$ 
and a global  action of $G$ on $Y,$ such that the initial partial action can be obtained as a restriction of the global one.\footnote{Such a  global action with a natural additional assumption  is called an enveloping action.} 
This was  studied first   in  
the PhD Thesis \cite{Abadie} (see also  \cite{AbadieTwo}) and, independently from   \cite{Abadie}, \cite{AbadieTwo}, in  \cite{KL} and \cite{St2}.   
Subsequent 
results    were obtained in \cite{Avi1},   \cite{BCFP},
 \cite{CF2},  \cite{DE}, \cite{DES2}, \cite{DdRS}, \cite{EGG}, \cite{F},   and more recently in    \cite{AbDoExSi},   \cite{BeCoFeFlo},   \cite{BeF}, \cite{CFMarcos},  \cite{Ferraro},
\cite{PiUz1} and \cite{PiUz2}.  The question  was also considered for partial semigroup actions in \cite{GouHol1}, \cite{Hol1}, \cite{Khry1}, \cite{Megre2}, \cite{Kud}, \cite{Limin},  
for partial groupoid actions in \cite{Gil, BP,BPi}, 
and around partial Hopf (co)actions in \cite{AAB,AB2,AB3,ABDP2,CasPaqQuaSant,CasPaqQuaSant2,CasQua}.\\

The importance of the  glo\-ba\-li\-za\-tion problem lies in the possibility to relate partial actions  with  global ones and this way try to move from global  results to the partial setting, producing  more general facts, as well as to obtain applications to  the global case in situations in which partial actions appear naturally, as it occurred in \cite{AraE1}.  Thus facts about globalization from    \cite{AbadieTwo} were used   in \cite{Norling2015}  with respect to  $K$-theory of  reduced $C^*$-algebras of  $0$-$F$-inverse semigroups,  in  \cite{LiNorling2016} in the $K$-theoretic study of reduced crossed products  attached to totally disconnected dynamical systems, and  in \cite{Artigue}  for  partial flows with application to  Lyapunov functions.   In addition, glo\-ba\-li\-za\-ble partial actions  were essential  for the development of Galois Theory of partial group actions in \cite{DFP}, 
for the elaboration of the concept of a partial Hopf (co)action in \cite{CaenJan},    as well as  in a series of  ring theoretic and 
Galois theoretic investigations  in \cite{Avi1}, \cite{AvFe}, \cite{AvFeLa}, \cite{AvLa},  \cite{AvLa3},   \cite{BLP}, \cite{BP},
\cite{BeCo}, \cite{BeCoFeFlo},  \cite{CaenDGr}, \cite{CCF}, \cite{C}, \cite{C2}, \cite{CF},  \cite{CFHM}, \cite{CFM},
\cite{FL},  \cite{FrP},  \cite{PRSantA}, \cite{PSantA}.\\

Another way to relate partial and global actions was given in  
R.~Exel's paper \cite{E1}, in which for any group $G$ a semigroup ${\mathcal S}(G)$ was defined by means of generators $\{ [g] \; | g\in G\}$ and relations:
\begin{eqnarray}
[g^{-1}][g][h]&=&[g^{-1}][gh],\nonumber\\
{}[g][h][h^{-1}]&=&[gh][h^{-1}],\nonumber\\
{}[g][1_G]&=&[g],\nonumber
\end{eqnarray}
$g, h \in G$ (it follows  that $[1_G][g]=[g]$).\footnote{The semigroup ${\mathcal S}(G)$ was denoted by $E(G)$ in \cite{DN}, \cite{DN2} and \cite{DoNoPi}.}  It was proved  that ${\mathcal S}(G)$ is an inverse semigroup \cite[Theorem 3.4]{E1},
 and     the partial actions of $G$ on a set $X$ are in one-to-one correspondence with the (global) actions of 
${\mathcal S}(G)$ on $X$ (see \cite[Theorem 4.2]{E1}). The injective map
$$G\ni g \mapsto [g] \in {\mathcal S}(G),$$ is the canonical  partial homomorphism, which  plays a key role in the above result.  Thus, instead of embedding $X$ into a larger set (or  object), one 
``expands''  $G$, obtaining a global action on  $X.$  It directly follows from the Exel's definition of ${\mathcal S}(G)$ that any partial representation of $G$ can be uniquely extended by means of the above mentioned 
map $G\to {\mathcal S}(G)$   to a (usual) representation of ${\mathcal S}(G)$ 
\cite[Proposition 2.2]{E1}.\\

The inverse semigroup $S(G)$ can be characterized as an {\it expansion} of $G$ in the sense of  J.-C. Birget, J. Rhodes  \cite{BR1} (see also \cite{BR2}), as follows. An expansion is defined as a  functor $F$ from the category of semigroups  into some special category of semigroups which has the property that there is a natural transformation $\eta $ from the functor $F$ to the identity functor such that $\eta _S$ is surjective for every semigroup $S.$  Amongst the several expansions  discussed in   \cite{BR1}, the so called {\it prefix expansion} is relevant for us.  In \cite[Proposition 1]{Sze} M. Szendrei gave a simple and very useful description of the prefix expansion ${\rm Pr}(G)$ (also denoted by $\tilde{G}^R$) of a group $G,$ and proved that 
${\rm Pr}(G)$ is an $F$-inverse semigroup which enjoys a certain universal ``$F$-inverse property'' 
\cite[Corollary 3]{Sze}. The latter is a consequence of a more general result \cite[Theorem 2]{Sze}, which states that $ {\rm Pr}(.)$ is a functor from the category of groups into the category of $F$-inverse
semigroups, which  is a left adjoint
of the functor assigning the greatest group homomorphic image to any $F$-inverse semigroup.\\

The construction can be applied to any monoid  $M$ (or even to any semigroup) and is called the  {\it Szendrei expansion} of $M$ \cite{FountainGomes},   \cite{GouHol1}, \cite{GouHol2}, \cite{Hol1}.   In general the Szendrei expansion ${\rm Sz}(M)$ of a monoid $M$ differs from the  Birget-Rhodes prefix expansion ${\rm Pr}(M),$ but  ${\rm Sz}(G)= {\rm Pr}(G)$  if $G$ is a  group. M. Szendrei's idea is of high  importance for partial actions, and, in oder to recall it,  denote by ${\mathcal P}_1 (M)$ the set of all finite subsets of a monoid $M$ containing $1.$ Then 
$${\rm Sz}(M) = \{(A, x) : A\in  {\mathcal P}_1 (M), x \in A \},$$ with the operation given by 
$$ (A, x) (B, y) = (A \cup  xB, xy).$$  If $M$ is  a semigroup without identity element,  then one adjoins an external $1$ and replaces $M$ by $M^1=M\cup \{1\}$ in the above definition.   J. Kellendonk and M. V. Lawson  have shown in   \cite{KL} that ${\mathcal S}(G)$ is isomorphic to ${\rm Sz}(G),$ and, as a consequence, to ${\rm Pr}(G).$ It is derived from  a universal property of 
${\rm Pr}(G)$  with respect to premorphisms \cite[Theorem 2.4]{KL}.  Thanks to the isomorphism 
${\mathcal S}(G) \cong  {\rm Pr}(G),$ Exel's definition of ${\mathcal S}(G)$ gives a presentation for   ${\rm Pr}(G)={\rm Sz}(G)$ in terms of generators and relations.\\

 The above  facts turned  
${\mathcal S}(G)$ into a highly important  tool, especially when dealing with partial projective group representations  \cite{DN}, \cite{DN2}, \cite{DoNoPi},  and the co\-ho\-mo\-lo\-gy theory based on partial actions  \cite{DKh},  \cite{DKh2},  \cite{DKh3}, as well as  when relating crossed products by partial actions of groups with crossed products by inverse semigroup actions \cite{EV}.    They were   used in \cite{ChLi}  to   study realizations of 
  ${\rm Pr}(G)$ as an inverse monoid of M\"obius type related to a partial action of $G$ on a Hausdorff space.  
Moreover, the expansion method was further developed and used for partial actions of inverse semigroups in \cite{LawMarSte}, \cite{BuE}, of  groupoids in 
\cite{Gil}, \cite{B},    \cite{BFP},     of monoids in \cite{Hol1}, of restriction semigroups in \cite{GouHol1} and    of inductive constellations in \cite{GouHol2}.  In particular, A.~Buss and R.~Exel gave in \cite{BuE} a presentation 
for the  generalized prefix expansion ${\rm Pr}(G)$ of an inverse semigroup $G,$ introduced earlier by    M. V.~Lawson, 
S. W. Margolis and B. Steinberg in \cite{LawMarSte}. It is proved in \cite{BuE} that twisted partial actions of $G$  on $C^*$-algebras correspond to twisted global actions of ${\rm Pr}(G),$ and  this correspondence preserves $C^*$-crossed products.\\

It became clear  already from the results in  \cite{Abadie}, \cite{AbadieTwo}, \cite{KL} and \cite{St2} that 
 the globalization problem strongly  depends on the category under con\-si\-de\-ra\-tion. In particular,  
 globalizations of
  partial actions on topological spaces   always exist, nevertheless, the to\-po\-lo\-gi\-cal properties of the initial space 
are not necessarily shared by the space under the global
  action. According to  \cite[Example 1.4]{AbadieTwo} there exists  a partial  group action on a
  Hausdorff space whose (minimal)  glo\-ba\-li\-za\-tion acts on a
  non-Hausdorff space, and, moreover,    in \cite[Proposition
  1.2]{AbadieTwo} a criteria was given for the preservation of the
  Hausdorff property under globalization. Because of the categorical
  equivalence between locally compact Hausdorff spaces and commutative
  $C^{\ast}$-algebras, this implies   that partial actions on
  $C^{\ast}$-algebras are not glo\-ba\-li\-za\-ble in ge\-ne\-ral (see
  Proposition 2.1 in \cite{AbadieTwo} for a criteria of the existence
  of a globalization of  a partial action on commutative
  $C^{\ast}$-algebras). On the other hand,  it was shown in
  \cite[Theorem 6.1]{AbadieTwo} that  
 they are globalizable ``up
  to Morita equivalence''. More precisely,  the concept of Morita equivalence of
  partial actions of locally compact groups on $C^{\ast}$-algebras was
  introduced and stu\-di\-ed in \cite{AbadieTwo} (see also \cite{Sieben2001}
  for the case of discrete groups), as well as that of a Morita enveloping action,  which is roughly a global action, whose restriction is a partial action Morita equivalent to the initial one.  It was shown that Morita equivalent partial actions have (strongly) Morita equivalent reduced $C^*$-crossed products.  Furthermore, the reduced $C^*$-crossed product of a partial action is (strongly) Morita equivalent to that of the Morita enveloping action. Notice that a particular  Morita equivalence fact based on a partial action on a commutative $C^*$-algebra from \cite{AbadieTwo} was recently related to a result from \cite{KoshkamSkandalis} by the authors of  the above mentioned paper \cite{MilanSt}. In the latter article a number of important Morita equivalence facts for $C^*$-algebras were estabilshed.\\

Influenced by F. Abadie's paper \cite{AbadieTwo}, the abstract ring theoretic analogues of the above mentioned concepts  from  \cite{AbadieTwo} were defined and 
studied in \cite{AbDoExSi}. Facts similar to those  from \cite{AbadieTwo} were proved  in the context of idempotent rings, whose Morita theory was developed in \cite{GS}. Moreover, some further Morita theoretic results were also obtained, including   the behavior of Morita equivalent partial actions  under the passage to matrices of infinite size with finite number of non-zero entries.  The latter  has no $C^{\ast}$-algebraic analogue so far, and  the treatment heavily depends on the technique worked out in \cite{DES1} to prove a ring theoretic  analogue of a stabilization result  for $C^*$-algebraic bundles from \cite{E0}.\\

The theory in \cite{AbDoExSi} is developed  for the so-called  regular partial group actions
on idempotent rings. This   includes all partial actions on  $C^{\ast}$-algebras, as well as all $s$-unital partial actions.  Note that the $s$-unital condition on a ring ge\-ne\-ra\-li\-zes all kind of unity conditions in ring theory, including   the 
existence of local units.
The regularity assumption imposes a mild restriction on the domains of  the partial isomorphisms involved in a partial action, more precisely,   the intersection of domains are assumed to coincide with their product.  This is a suitable constraint  since, 
on one hand, it resolves the discrepancy between the definitions of partial actions given   in \cite{DE} and \cite{DES1}, and on the other, in almost all  investigations on the subject the considered partial group actions on algebras (or rings) are regular, so that this concept  provides a sufficiently  general framework for the theory.\\

It is a   well-known fact that Morita equivalent commutative
rings with $1$ are necessarily isomorphic.  More generally this holds
for   non-degenerate idempotent rings, as   established in
\cite{GS}.  An  ana\-lo\-gous result for partial actions was given in  \cite{AbDoExSi}: Morita equivalent $s$-unital partial actions of a group $G$ on commutative algebras must be isomorphic.   A similar fact for $C^{\ast}$-algebras was also  established  in  \cite{AbDoExSi}: Morita equivalent partial actions of a discrete group $G$ on commutative $C^{\ast}$-algebras are necessarily isomorphic.\\

The theory of Morita equivalent partial actions on $C^*$-algebras from  \cite{AbadieTwo} was further developed by introducing the concept of a weak equivalence for arbitrary Fell bundles over  locally compact Hausdorff groups \cite{AbFerraro2} and that of their strong equivalence \cite{AbBussFerraro}. The weak equivalence captures the relation between a globalizable partial group action $\alpha$ on a $C^*$-algebra and its enveloping action $\beta$: the Fell bundles  corresponding to $\alpha$ and $\beta $ are weakly equivalent (see \cite[Example 2.21]{AbFerraro2}). More generally, it follows from the results in    \cite{AbadieTwo} and  \cite{AbFerraro2} that every Fell bundle associated to a partial action is weakly equivalent to the Fell bundle associated to a global action \cite[Corollary 5.15]{AbFerraro2}. Furthermore, 
amenability is preserved under weak equivalence   \cite{AbFerraro2} (see also  \cite{AbBussFerraro}).\\

The more restrictive notion of a strong equivalence of Fell bundles  is a natural generalization of that of a Morita equivalence of partial actions:   two partial actions on $C^*$-algebras are Morita equivalent if and only if the Fell bundles associated to them are  strongly  equivalent  \cite[Corollary 4.9]{AbBussFerraro}.
 One of the main results in \cite{AbBussFerraro} asserts that every Fell bundle is strongly  equivalent to a semidirect product Fell bundle for a partial  action. Consequently, every Fell bundle is weakly equivalent to the semidirect product Fell bundle of a global action. In addition, it is shown in  \cite[Corollary 4.3]{AbBussFerraro} (see also \cite[Proposition 4.13]{AbFerraro2}) that  weakly equivalent Fell bundles have (strongly) Morita equivalent full and reduced cross-sectional $C^*$-algebras. Notice that \cite[Proposition 7.1]{KwaMeyer} implies that the reduced $C^*$-algebra of a Fell bundle over a discrete group $G$ is Morita equivalent to the reduced crossed product by a global action of $G$.\footnote{The author thanks Fernando Abadie for drawing his attention to this fact.}\\

In the recent article  \cite{Ferraro} the problem of the existence of  an en\-ve\-lo\-ping action for  a partial group action  on  a non-necessarily unital 
(abstract) ring was investigated and applied to the globalization problem for partial actions on  $C^*$-algebras and equivalence  Hilbert bimodules.  More spe\-ci\-fi\-cal\-ly, amongst various facts it was proved  that   if a partial 
action $\alpha = \{ \alpha _g : A_{g\m} \to A_g, g\in G\}$ of a (discrete) group on a (non-necessarily unital) ring $A$ admits a globalization, then 
the following condition is satisfied:\\

\noindent For each $(g,a,b) \in  G \times A \times A$    there exists   $u \in A_g$  such that 
\begin{equation}\label{FerraroCondition}
  cu = \alpha _g (\alpha _{g\m} (c)a)b \;\text{ and }\; uc = \alpha _g (a\alpha_{g\m} (bc)),
\end{equation}  for all  $c \in A_{g}$ (see   \cite[Theorem 2]{Ferraro}).\\

\noindent Under  appropriate non-degeneracy assumptions (which always hold for $C^*$-algebras) condition (\ref{FerraroCondition}) becomes sufficient for the existence of a globalization for $\alpha .$\\

 For the $C^*$ case, let $G$ be  an arbitrary topological  group  and $\alpha $ be a $C^*$-partial action of $G$ on a $C^*$-algebra $A.$  Contrary to the usual practice in $C^*$-theory, $G$ is not assumed to be  Hausdorff nor locally compact. The topological generality of $G$ aims to understand the effect of the group's topology on the existence of a globalization.  Write $G^{\rm dis}$ regarding $G$ as a discrete group and let $\alpha ^{\rm dis}$ be the 
 $C^*$-partial action of  $G^{\rm dis}$ on $A$ given by $\alpha .$  Then \cite[Corollary 2]{Ferraro} means that for the globalization problem the topology of $G$ can be forgotten, i.e. $\alpha $ admits a $C^*$-globalization if and only of so does   $\alpha ^{\rm dis}.$ Moreover, \cite[Theorem 5]{Ferraro} asserts 
 that $\alpha $ has a $C^*$-globalizaton if and only of $\alpha $ admits a ring theoretic globalization. The latter is shown to be equivalent to the condition given by the first equality in (\ref{FerraroCondition}). As a consequence, it follows using \cite[Theorem 4.5]{DE} that a $C^*$-partial action $\alpha $ on 
 a unital $C^*$-algebra $A$ possesses a $C^*$-globalization if and only if $\alpha $ is unital, spelling out thus  a $C^*$-version of \cite[Theorem 4.5]{DE}.\\
 
 The above mentioned  Morita theory  for partial actions on $C^*$-algebras  \cite{AbadieTwo} involves the so-called $C^*$-ternary rings, which also may be  seen as equivalence Hilbert bimodules, and the author  of \cite{Ferraro} proves that  a partial action of a topological group on an equivalence Hilbert bimodule has a glo\-ba\-li\-za\-tion if and only if its linking partial action   \cite{AbadieTwo} has a $C^*$-globalization \cite[Corollary 6]{Ferraro}. Notice that new applications of $C^*$-ternary rings to $C^*$-algebras were given in \cite{AbadieFerraro}.\\

Another recent development on the globalization problem was obtained in \cite{ABDP2} with respect to twisted Hopf partial actions. Partial actions of Hopf algebras on algebras were defined by S.~Caenepeel and K.~Janssen in \cite{CaenJan} influenced 
by the Galois theory of partial group actions  developed in \cite{DFP}.  The latter stimulated also   further Galois theoretic results in \cite{CaenDGr}, which were based on a coring ${\mathcal C} $ constructed for a unital partial action of a 
finite group, offering thus a more conceptual approach to partial Galois theory via Galois corings.\footnote{The authors of  \cite{CaenDGr}  use the term  idempotent partial action, nevertheless we prefer to employ the latter name in a different sense.}   The coring ${\mathcal C} $ was shown to fit the general theory of cleft bi\-co\-mo\-du\-les  in \cite{bohmverc}, and, in addition, in \cite{Brz}  descent theory for corings was  applied, using   ${\mathcal C} ,$ to define  non-Abelian Galois cohomology  ($i=0, 1$) for unital partial Galois actions of finite groups. The article by  S.~Caenepeel and K.~Janssen  \cite{CaenJan}, in its turn,  became  the starting point for a series of  investigation of partial Hopf (co)actions (see the articles cited above), in particular, several globalization results were obtained in  \cite{AAB}, \cite{AB2}, \cite{AB3}, \cite{CasPaqQuaSant}, \cite{CasPaqQuaSant2}, \cite{CasQua} and     \cite{KrakenQuastTimm}.\\

On the other hand, the twisted version of partial group actions on (abstract) algebras were introduced and studied in \cite{DES1}. This was inspired by the R.~Exel's  notion of a continuous twisted partial action of a locally compact group on a
$C^*$-algebra (a twisted partial $C^*$-dynamical system) and that of the cor\-res\-pon\-ding crossed pro\-duct \cite{E0}. 
  The new construction permitted one to show that any se\-cond countable $C^*$-algebraic bundle, which satisfies a certain regularity condition (automatically verified  if the unit fiber algebra is stable), is a $C^*$-crossed product of the unit fiber algebra by a continuous partial action of the base group  \cite{E0}. The algebraic version of the latter fact was established in \cite{DES1}.  This algebraic concept  was applied to Hecke algebras in \cite{E3}, where, among other results, it was proved that given a field $\kappa $
of characteristic $0,$ a group $G$ and subgroups $H, N \subseteq G$ with $N$  normal in
$G$ and $H$  normal in $N,$ there is a twisted partial action $\theta $ of $G/N$ on the group
algebra $\kappa (N/H)$ such that the Hecke algebra ${\mathcal H}(G,H)$ is isomorphic to the
crossed product $\kappa (N/H) \ast _{\theta} G/N.$ The globalization problem for twisted partial group actions was investigated 
in \cite{DES2}, whereas other algebraic results 
on twisted partial actions on rings and corresponding crossed products were obtained in \cite{BLP}, \cite{BaraCoSo},   \cite{BeCo}, \cite{BeCoFeFlo} and \cite{PSantA}.\\

Motivated by the concept of a twisted partial group action given in \cite{DES1} on one hand, and the developments on partial Hopf actions on the other, twisted partial actions of Hopf algebras on rings were introduced in \cite{ABDP}, as well as the corresponding crossed products. Examples using  algebraic groups were elaborated, more precisely,  
actions of an affine algebraic group on affine varieties give rise to  coactions of the corresponding commutative Hopf algebra $H$  on  the coordinate algebras of the varieties,  restrictions of which produce  concrete examples of partial Hopf coactions.  Then one  dualizes in order to obtain  partial Hopf actions.  The dualization passage works theoretically, but may far from being easy in concrete examples.  One possibility is to try to identify the finite dual $H^0$ for a specific $H$ obtained in the above 
way.  A more flexible possibility is to find a concrete Hopf algebra $H_1$ such that $H$ and $H_1$ form a dual pairing.  Then  \cite[Proposition 8]{AB2} produces a partial action of $H_1.$ A concrete example was elaborated  \cite{ABDP}, which was easily endowed with a twisted structure by means of a $2$-cocycle. Furthermore,  symmetric 
twisted partial Hopf actions were introduced in   \cite{ABDP} in order to treat the convolution invertibility of the partial cocycle in a ma\-na\-ge\-a\-ble way. This permitted one to get a relation between crossed products and the so-called partially cleft extensions. The latter were  also introduced in the same paper    \cite{ABDP},  and the definition  reflects  ``partiality'' in more than one way, incorporating, in particular, some equalities 
already proved to be significant in the study of partial group actions and partial representations. Then one of the main   
facts  in     \cite{ABDP} states that the partial cleft extensions over the coinvariants $A$ are exactly the crossed pro\-ducts by symmetric twisted partial Hopf actions on $A.$ One should notice that it was proved in \cite{FerGonRod} that the crossed pro\-ducts by twisted partial Hopf actions \cite{ABDP} form a particular case of the more general weak crossed products defined in \cite{FerGonRod0}.\\

 Then it became natural to investigate  the globalization problem for twis\-ted partial Hopf actions on rings \cite{ABDP2}. The main result says that  a symmetric twisted partial action of a Hopf algebra $H$ on a unital algebra $\A$ associated to the symmetric pair of partial cocycles $\omega$ and $\omega '$,  is globalizable if, and only if, there exists a normalized convolution invertible linear map $\tilde{\omega} :H\otimes H\rightarrow \A$ satisfying certain compatibility conditions, intertwining the partial action of $H$ on $\A$ and the restriction of the twisted action of $H$ on $\B$. A series 
of examples were elaborated. In particular, a complete characterization of partial measuring maps were given 
in the case of $H$ being the group algebra $\kappa G$, the dual $(\kappa G)^*$ of the group algebra  $ \kappa G $ of a 
finite group $G,$ and the Sweedler Hopf algebra $H_4,$ and $\A$ coinciding (in all three cases) with  the base field $\kappa $.
Furthermore, symmetric twisted partial Hopf actions were  described in details for specific Hopf algebras.
 In particular,  the case of a group algebra $\kappa G$ recovers  the theory of twisted partial actions of groups as developed 
in \cite{DES1} and \cite{DES2}. For the Sweedler Hopf algebra $H_4$, the only symmetric twisted partial actions are the global ones. In addition, for the case of the dual  $(\kappa G)^*$ of the group algebra of a finite group $G$, the partial 
 appearing cocycles  have remarkable symmetries, and they were related  to  global cocycles of  dual group algebras of quotient groups. The specific example for the Klein four-group $K_4$, acting on the base field $\kappa ,$ was done in more details, 
and it was  shown that 
the symmetric twisted partial actions of $(\kappa K_4)^{\ast}$ on $\kappa $ are parametrized by the zeros
 $(x, y) \in \kappa ^2$ of a polynomial in $x, y$ of degree $2.$ An explicit partial cocycle for   $(\kappa K_4)^{\ast}$ was given which leads to a globalizable symmetric 
twisted partial action. Moreover, it was shown that  the example of a twisted partial Hopf action, constructed in \cite{ABDP} using the relation between algebraic groups and Hopf algebras, is globalizable.\\ 

The algebraic advances on Hopf partial actions influenced a $C^*$-theoretic development: in the already mentioned paper   \cite{KrakenQuastTimm}, the approach for partial Hopf coactions from \cite{AB3} motivated a globalization result  for partial coactions of $C^*$-quantum groups satisfying a mild restriction, which always holds if the quantum group is discrete, or if the $C^*$-algebra of the quantum group is nuclear. The construction of the enveloping action  gives  a left adjoint to the forgetful functor from coactions to partial coactions (with an appropriate choice of categories).\\

The algebraic concept of twisted partial actions also motivated the study of projective partial group representations,  the corresponding partial Schur Multiplier and the relation  to partial group actions with $\kappa $-valued twistings in \cite{DN}, \cite{DN2} and \cite{DoNoPi}, and contributing thus towards the elaboration of a background for a general cohomology theory based on partial actions.\\

  A (usual) projective representation of a group $G$ can be defined as  a homomorphism from $G$ to the projective linear group ${\rm PGL}_n(\kappa ).$ In order to define  partial projective representations one replaces usual homomorphisms by partial ones as follows.  Denote by ${\rm PMat}_n \kappa $  the  monoid of the projective   $n\times n$ matrices over a field $\kappa $, i.e. 
${\rm PMat}_n \kappa = ( {\rm Mat}_n \kappa ) /\lambda ,$ where  $\lambda $ is the congruence given by   
$A\lambda B  \Longleftrightarrow $  $A=cB$ for some  $c\in \kappa ^{\ast}$. Then we define a {\it partial projective representation} of $G$ as a partial homomorphism of the form $ G \to {\rm PMat}_n \kappa$ \cite{DN}. Taking representatives in the congruence classes, we may consider a partial projective representation as a function of the form  
$ G \to {\rm Mat}_n \kappa .$
As in the classical case, factor sets appear naturally: a factor set of a partial projective representation $\Gamma : G \to   {\rm Mat}_n \kappa $ is a function 
$\sigma : G \times G \to \kappa $ such that 
\begin{equation}\label{parfac1}\sigma (g,h) =0 \Longleftrightarrow \Gamma (g) \Gamma (h) =0,\end{equation}  
\begin{equation}\label{parfac2}\Gamma (g\m ) \Gamma (g) \Gamma (h) =  \Gamma (g\m ) \Gamma (g h) \sigma (g,h),\end{equation} and
\begin{equation}\label{parfac3}\Gamma (g ) \Gamma (h) \Gamma (h\m ) =  \Gamma (g h ) \Gamma ( h\m) \sigma(g,h), \end{equation}    for all $g,h \in G$ (see  \cite[Theorem 3]{DN}).
The set $$X=X_{\sigma} = \{ (g,h) \in G\times G \, : \, \Gamma (g) \Gamma (h) \neq 0 \}$$ is called the {\it domain} of $\sigma $
(or the domain of $\Gamma$).\\

More generally, in \cite{DN} the partial projective representations are con\-si\-de\-red as functions of the form 
$\Gamma : G \to M,$ where $M$ is a so-called $\kappa $-cancellative monoid (see \cite[Definition 2]{DN}).\\

The theory of partial projective representations is strongly related to  Exel's semigroup ${\mathcal S}(G).$ In fact, they can be 
alternatively defined via projective representations of ${\mathcal S}(G),$ so that the theory  of projective representations of semigroups and their Schur multipliers,  elaborated by B. Novikov in  \cite{nov1},  \cite{nov}, \cite{nov3} (see also  \cite{nov4}),  comes into the picture as an essential working tool. The usual cohomology of semigroups does not serve the projective semigroup representations, instead the more general $0$-cohomology  \cite{nov3} fits them with its natural partial flavor.\\ 

The factor sets  of the partial projective representations of $G$ form 
   a commutative  semigroup, whose equivalence classes
constitute the  partial Schur Multiplier $pM(G)$. The latter is a commutative inverse semigroup, and as such it  is a semilattice  of abelian groups $pM_X(G),$ called components, where $X$ runs over the domains of the partial projective representations of $G.$  One of the component, namely, the group of the equivalence classes of the totally defined factor sets $pM_{G\times G}(G),$  contains the usual Schur Multiplier  $M(G)$ of $G$, but  $pM_{G\times G}(G)$ is essentially bigger than $M(G).$ In \cite{DN2} the structure  $pM_{G\times G}(G)$ was 
investigated over an algebraically closed field $\kappa $ and the technique of  \cite{DN2} was further developed in  \cite{DoNoPi}, permitting one to extend the description to any component  $pM_X(G),$  and to show that each   $pM_X(G)$ is an epimorphic image of  $pM_{G\times G}(G).$ Furthermore, it was also shown in \cite{DoNoPi} that each 
  $pM_X(G)$ is an epimorphic image of a direct power of $\kappa ^{\ast}$ ($\kappa $ is assumed to be algebraically closed).\\  

Since the partial projective  representations of a group $G$ are intimately related to the projective representations of the inverse semigroup ${\mathcal S}(G),$ it was natural to explore this connection and use it as a working tool.
  This relation was established in \cite{DN} and further  explored in  \cite{DN2} and  \cite{DoNoPi}, incorporating some simplifications. One of them  is  a passage from the semigroup ${\mathcal S}(G)$ to its quotient ${\mathcal S}_3(G).$  In order to define the  latter, we identify ${\mathcal S}(G)$ with $ {\rm Sz}(G),$   sending  the generator $[g]\in {\mathcal S}(G), g\in G,$ from the initial definition of ${\mathcal S}(G)$  to the pair $(\{1,g\},g) \in {\rm Sz}(G),$ \cite{KL}. Denoting by  $N_k$ the ideal $\{(R,g)\in {\mathcal S}(G)\ |\ |R|\geq k+1\}$ of  ${\mathcal S}(G),$ the semigroup 
 ${\mathcal S}_3(G)$ is the quotient ${\mathcal S}(G)/N_3(G).$\footnote{The semigroup  ${\mathcal S}_3(G)$
  was denoted by $E_3(G)$ in  \cite{DN2} and \cite{DoNoPi}.} Reductions based on this passage   permitted one to obtain  characterizations of  partial factor sets of $G$ over an algebraically closed field, one of which gives defining equalities for the partial factor sets that  curiously  have nothing to do with the $2$-cocycle identity, but nicely   incorporate  the symmetry under the action of the symmetric group $S_3$  
\cite[Theorem 5.6]{DoNoPi}.  These results  led to  general facts  on the structure of the components of    $pM(G)$  and allowed one to perform calculations for concrete groups.\\

R.~Exel's  concept  of a continuous twisted partial action \cite{E0} and its purely ring theoretic version \cite{DES1} 
involve a general twisting which satisfies the $2$-cocycle identity in some restricted sense, and it was natural to fit this in some cohomology theory. A relation of this cohomology to the partial Schur Multiplier is naturally expected, as the classical 
 Schur Multiplier   is isomorphic  to the cohomology group $H^2(G, \bbC ),$ where the action of $G$ on $\bbC$
 is trivial. The main  idea is to replace  usual $G$-modules, i.e. global actions of $G$ on abelian groups, by partial 
$G$-modules, which are partial actions of $G$ on commutative monoids. The first step was done in \cite{DKh} with more recent developments in the preprints \cite{DKh2} and \cite{DKh3}.  Since  the twistings  in \cite{DES1}  take  values in  multiplier algebras of products of some ideals, it was reasonable 
to avoid  multipliers at the beginning, imposing a rather usual  restriction on the partial action, namely, that it is unital. We recall that a partial action on a ring (or a semigroup) is called unital if each domain is an ideal generated by an idempotent, which is central in the ring (semigroup).\\ 

 If we assume that the ring  $\A$ under a unital  twisted partial action  of a group $G$ is  commutative, then the action 
  falls into two parts: a partial action $\alpha $ of $G$ on 
$\A$ and its twisting. This way we derive the notion of  a partial $2$-cocycle (the twisting) whose  values belong to groups of invertible elements of appropriate ideals of $\A.$ The concept of a partial $2$-coboundary then follows from that of an equivalence of twisted partial actions introduced in~\cite{DES2}.  Of course, in the general definition we do not need a ring structure on $\A$, so assuming that $A$  is a commutative multiplicative monoid, one comes to the definition of the second co\-ho\-mo\-lo\-gy 
group $H^2(G,A).$  The groups $H^n(G,A)$ with arbitrary $n$ are defined in a similar way. More precisely, given a partial $G$-module $A,$ i.e. a commutative monoid endowed with a partial action $\alpha =\{ \alpha_g : A_{g\m} \to A_g, g \in G\},$ 
 for any $n>0$ one defines  {\it $n$-cochains} of $G$ with  values in $A$ as functions $f:G^n\to A$, such that $f(x_1,\dots,x_n)$ is an invertible element of the ideal $A_{(x_1,\dots,x_n)}=A_{x_1}A_{x_1x_2}\dots A_{x_1\dots x_n}$. 
By {\it a $0$-cochain} we  mean an invertible element of $A$. Then the set  $C^n(G,A)$ of $n$-cochains is an abelian group under the pointwise multiplication with the identity  element
\[
	e_n(x_1,\dots,x_n)=1_{x_1}1_{x_1x_2}\dots 1_{x_1\dots x_n},
\]
and the inverse of $f\in C^n(G,A)$ being $f^{-1}(x_1,\dots,x_n)=f(x_1,\dots,x_n)^{-1}$, where $f(x_1,\dots,x_n)^{-1}$ means the inverse of $f(x_1,\dots,x_n)$ in $A_{(x_1,\dots,x_n)}$.\\
  
Next, for any $f\in C^n(G,A)$ and $x_1,\dots,x_{n+1}\in G$ define the {\it  coboundary map}: 
	\begin{align*}\label{eq-coboundary_hom}
		(\delta^nf)(x_1,\dots,x_{n+1})&=\alpha_{x_1}(1_{x_1^{-1}}f(x_2,\dots,x_{n+1}))\notag\\
&\prod_{i=1}^nf(x_1,\dots , x_ix_{i+1}, \dots,x_{n+1})^{(-1)^i}\notag\\
&f(x_1,\dots,x_n)^{(-1)^{n+1}}. 
	\end{align*}
Here the inverse elements are taken in the corresponding ideals. If $n=0$ and $a$ is an invertible element of $A$, we set $(\delta^0a)(x)=\alpha_x(1_{x^{-1}}a)a^{-1}$. Then  $\delta^n$ is a homomorphism $C^n(G,A)\to C^{n+1}(G,A)$, and 
	\begin{equation*}\label{eq-deltasquare}
		\delta^{n+1}\delta^nf=e_{n+2}
	\end{equation*} 
for any $f\in C^n(G,A)$ (see \cite{DKh}). Now,  as in the classical case, we define the abelian groups $Z^n(G,A)={\rm Ker}\, {\delta^n}$, $B^n(G,A)={\rm Im}\, {\delta^{n-1}}$ and $$H^n(G,A)={\rm Ker}\, {\delta^n}/{\rm Im}\, {\delta^{n-1}}$$ of {\it partial $n$-cocycles, $n$-co\-boun\-da\-ries and $n$-cohomologies} of $G$ with values in $A$, $n\ge 1$ $(H^0(G,A)=Z^0(G,A)={\rm Ker}\, {\delta^0})$.\\

One actually may replace $A$ by an appropriate  submonoid $\tilde{A}$, which is inverse \cite{DKh}. This brings $A$ closer to the classical case, as the commutative 
inverse monoids are natural generalizations of abelian groups, not being too far from them.\\

One of the difficulties with the partial $G$-modules is that they do not form an abelian category. Nevertheless, our partial 
cohomology   can be related to the Lausch-Leech-Loganathan cohomology of inverse semigroups (see~\cite{Lausch}, \cite{Leech} and~\cite{Loganathan}) via the R.~Exel's inverse monoid ${\mathcal S}(G).$ From a unital partial action of $G$ on $A$ one comes to an action of ${\mathcal S}(G)$ and then to an ``almost'' Lausch's ${\mathcal S}(G)$-module structure on $A$. The latter can be seen as a mo\-du\-le in the sense of H.~Lausch over an epimorphic image of ${\mathcal S}(G)$, provided that $A$ is an inverse partial $G$-module. Thus our category is made up of abelian ``pieces'' which are categories of Lausch's modules over epimorphic images of ${\mathcal S}(G)$. This way we are able to define free objects and free resolutions which lead to $H^n(G,A)$ \cite{DKh}. We also showed that the partial Schur multiplier  $pM(G)$ is a union of $2$-cohomology groups of $G$ with values in non-necessarily trivial partial $G$-modules.\\

Then it was natural  to give an interpretation of  the  partial $2$-cohomology group in terms of extensions. 
This was initiated  in \cite{DKh2}, and  at some point it became clear that it is preferable to abandon the restriction on a partial 
action to be unital, imposed by cohomology theory in \cite{DKh}, covering thus a more general situation.    The key notion  is that  of an extension 
of a semilattice of groups $A$ by a group $G$ \cite{DKh}, the main   example being  the crossed product   $A\ast _\alpha G$ by a twisted partial action $\alpha $ of $G$ on $A$. Such an extension is  related to the notion of an extension 
of $A$ by an inverse semigroup $S,$  that of a twisted $S$-module and the corresponding crossed product given in \cite{Lausch}.   In particular, given an extension
$A \to U \to G$ there is a refinement 
$A \to U \to  S \to G$ such that $A\to U \to S$ is an extension of $A$ by $S,$ where $S$ is an $E$-unitary semigroup.
In order to make the theory work well we impose an admissibility condition on  the extensions   $A \to U \to G,$ and the 
$S$-module structures  on $A,$ that we obtain this way, possesses twistings which satisfy a normality condition, considered by 
N. Sieben in   \cite{Sieben98}, which is stronger than the one imposed by H. Lausch in \cite{Lausch}. For this reason we call them Sieben's twisted modules. Then we are able to show that any admissible extension $A \to U \to G$ is equi\-va\-lent to some crossed product  extension  $A \to A*_\alpha G \to G,$  and  the final fact in \cite{DKh2}  establishes  an equi\-va\-len\-ce preserving  one-to-one correspondence between twisted partial actions of groups on $A$ and Sieben's twisted module structures  on $A$ over $E$-unitary inverse semigroups.\\

 The theory of  extensions in \cite{DKh2} went beyond the cohomology theory in \cite{DKh}, so that it became clear that it is more natural to deal with a more general  multiplier 
valued cohomology theory. This is done  in \cite{DKh3} where  the cohomology groups    $H^n(G,A),$ with $A$ 
being a non-necessarily unital  partial $G$-module, were defined. For an inverse semigroup $S$ and an $S$-module structure on $A$ we also define the cohomology groups $H^n_\le(S^1,A^1)$ based on order preserving cochains and relate them 
to     $H^n(G,A).$  This is motivated by the fact that Sieben's twisted $S$-modules have order preserving twistings.
The elements of  the second cohomology group  $H^2(G,A)$ are proved to be in one-to-one correspondence  with the equivalence classes of the extensions of $A$ by $G.$ In addition, we define the concept of a split extension $A \to U \to G$
and prove that the  elements of    $H^1(G,A)$ are  in one-to-one correspondence  with the equivalence classes of splittings of $U.$\\

Note that partial group cohomology turned out to be  useful to study ideals of global  reduced $C^*$-crossed products:  to a given global $C^*$-dynamical system  the authors of  \cite{KennedySchafhauser}
associate a partial $C^*$-dynamical system,  giving rise to 
a ``twisted partial representation'', which is a projective partial representation, whose factor set (twist) is a partial $2$-cocycle $\sigma$. If $\sigma $  is a $2$-coboundary, then the initial (global)   $C^*$-dynamical system is said to have ``vanishing obstruction''. As the authors say: ``In a certain precise sense, the twist is the only obstruction to understanding the ideal structure of the reduced crossed product''.  Assuming  the vanishing obstruction property, several  necessary and sufficient conditions for the ideal intersection property are given.  The latter  is known to be closely related to the ideal separation property.\\

 Further results on partial projective representations were obtained in \cite{DdLP},   \cite{DoNoPi}, 
\cite{LimPin}, \cite{NovP},  \cite{Pi}, \cite{Pi2},  \cite{Pi4} and \cite{Pi5} (see also the short survey \cite{Pi3}). In particular, computations of the  partial Schur multiplier of concrete  groups in  \cite{DoNoPi},  \cite{LimPin}, \cite{NovP}, \cite{Pi2},    \cite{Pi4} show that each component is, in fact,  isomorphic to a  direct power of $\kappa ^{\ast},$   suggesting that  this should be true for all  groups. This motivated  the recent preprint \cite{DoSa}, in which this conjecture  was confirmed for all finite groups over an algebraically closed field. This surprisingly gives  a better understanding of the structure of   $pM(G)$ than one has for that of the usual  Schur Multiplier.\\

 Using \cite[Proposition 2]{DN} one can show that a partial factor set $\sigma $  satisfies the following weak $2$-cocycle condition:
 \begin{equation}\label{weak}
1\in\{x,y,z,xy,yz,xyz\}\Rightarrow\delta^2\sigma(x,y,z)=0.
\end{equation} The key idea in  \cite{DoSa}  is to replace  $\kappa ^{\ast}$ by an arbitrary abelian group $A$ and define {\it pre-cocycles}  which are functions $\sigma : G \times G \to A$ obeying condition (\ref{weak}).  They form a group denoted by $pZ^2(G,A).$ Since the equivalence of partial factor sets is defined modulo classical co\-boun\-da\-ri\-es, it is reasonable to introduce the factor group   $pZ^2(G,A)/B^2(G,A),$ which is denoted by $pH^2(G,A)$ and called the {\it pre-cohomology group} of $G$ with values in $A.$ Aiming to ge\-ne\-ra\-li\-ze 
the components   $pM_X(G)$ of the partial  Schur  multiplier  $pM(G),$ one wishes to restrict  pre-cocycles to domains. 
By \cite[Proposition 5.3]{DoNoPi} there is a bijection between the domains of factor sets of partial projective  representations of $G$ and the following set of ideals: 

\begin{equation*}
  \Lambda=\{I \trianglelefteq {\mathcal S}(G) \ |\ N_3\leq I \neq  {\mathcal S}(G)\}. 
\end{equation*} The correspondence takes a domain $X$ to the ideal $I\in \Lambda $ such that for all $x,y \in G:$
$$(x,y) \in X \Longleftrightarrow [x][y] \notin I. $$  

Let now $A\cup\{0\}$ be the semigroup obtained by adjoining a zero $0$ to the group $A$ (we assume the multiplicative notation for $A$).\footnote{In   \cite{DoSa} $A$ is additive and the zero element is denoted by $\infty.$}
Considering for any  ideal $I \in \Lambda$ the function  $\varepsilon_I:G\times G\to A\cup \{0\}$ determined by
 \[\epsilon_I(g,h)=\left\{\begin{array}{cl} 0 & \text{if} \;\;\; [g] [h]\in I,\\1_A&\mbox{otherwise,}\end{array}\right.\]
and setting 
\[Z^2(G,I;A)=pZ^2(G,A) \varepsilon_I,\  \ \ B^2(G,I;A)=B^2(G,A) \varepsilon_I, \] the quotient \[H^2(G,I;A)=Z^2(G,I;A)/B^2(G,I;A),  \]
is called the \emph{second partial cohomology group relative to} $I$ with coefficients in $A$. As the  abelian groups  $pM_X(G)$ assembles into the partial Schur  multiplier  $pM(G),$ it is natural to define 
the \emph{second partial cohomology semilattice of groups of} $G$ by
\[H^2(G,\Lambda;A)=\coprod_{I\in\Lambda}H^2(G,I;A).\] 

Another fruitful   point  in  \cite{DoSa} is to view partial factor sets as liftings (sections) $\varphi : G \to E$ of partial homomorphisms $\psi : G \to M$ related to appropriate central extensions $\kappa ^{\ast} \to E \to M$ fitting the following   commutative diagram:   

\begin{equation}\label{Eq:diagram}
\xymatrix{
&&G\ar[d]^\psi\ar@{.>}[dl]_\varphi\\
\kappa ^{\ast} \ar[r]^\iota&E\ar[r]^\pi&M
}\end{equation}  Then for any $\sigma \in  Z^2(G,I;A)$ it is possible to construct a central extension $A \to  E \to M,$ a partial  homomorphism $\psi : G \to M$   and its lifting (section) $\varphi :G \to E,$ such that the diagram  \eqref{Eq:diagram} is commutative upon the replacement of $\kappa ^{\ast}$ by $A,$ and $\sigma $ plays the role of the  factor set related to $\varphi $ as in (\ref{parfac1}),  (\ref{parfac2}), (\ref{parfac3}). This leads to a bijection  between the elements of $H^2(G,I;A)$ and the equivalence classes of certain appropriately defined central extensions \cite[Theorem 3.3]{DoSa}. This way the partial Schur multiplier   $pM(G)$ over an arbitrary field $\kappa $ coincides with  the second  partial cohomology semilattice of groups with values in  $A= \kappa ^{\ast},$ the components of which being second   relative partial cohomology groups.\\

Theorem 4.3 em \cite{DoSa} gives a  precise structure of the pre-cohomology group for a finite group $G$ of order $n$ an any abelian group $A:$

 \[pH^2(G,A)\simeq A^{m-n}\oplus A\otimes G/[G,G], \]
 where $m={(n^2+2|G_{(3)}|+3|G_{(2)}|+5)}/{6}$ and $G_{(k)}$ stands for the set of elements of order $k$ in $G$. 
Notice that $Z^2(G,N_3;A)=pZ^2(G,A).$ The structure result \cite[Theorem 5.2]{DoSa}  for a general relative cohomology group $H^2(G,I;A)$ is less precise but it is good enough to conclude that if $G$ is finite, $A= \kappa ^{\ast}$ and $\kappa $ is algebraically closed, then each 
 $H^2(G,I; \kappa ^{\ast})$ is isomorphic to a finite direct power of $ \kappa ^{\ast}.$\\

The relative partial cohomology  is consistent with the  partial co\-ho\-mo\-lo\-gy  in  \cite{DKh}: each second partial cohomology group relative to an ideal is isomorphic to the second partial cohomology group with values in an appropriate unital partial $G$-module \cite[Theorem 6.1]{DoSa}. As a consequence, one concludes that for a finite $G$ and an arbitrary field $\kappa ,$  each component of  $pM(G)$ is isomorphic to a single  cohomology group with values in some partial $G$-module. The latter  fact  brings  the partial Schur theory closer to the classical one.\\

While  partial group actions and partial group re\-pre\-sen\-ta\-ti\-ons are strongly related to inverse semigroups,
a generalization of the latter, namely, the (two-sided) restriction semigroups, found interesting applications of partial actions of monoids  in 
\cite{CornGould} and \cite{Kud}.  The notion of a (two-sided) restriction semigroup  is based  on the existence of  two unary operations $u \mapsto u^+$ and $u \mapsto u^*$ , which resemble  the maps $x \mapsto x\m x$ and $x \mapsto x x\m$ in an inverse semigroup. Its axiomatic definition incorporates the defining identities  of both left and right restriction semigroups, as well as two connecting relations $(u^+)^* = u^+,$  $(u^*)^+ = u^*$ (see \cite{CornGould} or \cite{Kud}).  In particular,  every subsemigroup of an inverse semigroup that is closed under $^+$  and $^*$ serves as an example of a  restriction semigroup. 
The  role of the 
relevant class of the $E$-unitary inverse semigroups is played now by the proper  restriction semigroups, for which in
 \cite{CornGould} a structure theorem  in terms of double partial actions of monoids on semilattices was given. This extends classical results  for inverse and ample semigroups. The double partial action from \cite{CornGould} can be reformulated in terms of only one partial action, and  
in \cite{Kud} classes of proper restriction semigroups determined by the properties of this partial action where classified.
In particular, a new important class of proper restriction semigroups was  introduced this way, which were called ultra proper restriction semigroups, and which nicely relates to  other well-established classes. The author of  \cite{Kud} uses partial actions to obtain various relevant results, in particular,    globalizations of a partial actions are applied  to establish a
McAlister-type theorem, as well an embedding fact into  $W$-products, extending  known results, and producing new and simpler proofs.\\

Another recent application of partial actions to semigroup theory was given in \cite{Khry1}, offering  a simple proof 
of a weakened version of  an  embedding theorem by O'Carroll \cite{O'Carroll}. The notion of a partial crossed product semigroup (also called partial  semidirect product, if the twisting is tri\-vi\-al),  considered earlier with some variations\footnote{In the case of a partial action on a semilattice, the  corresponding semidirect product construction was  essentially known much earlier without using the notion of a partial action (see, in particular, \cite{O'Carroll},   \cite{PetrichReilly}).} in \cite{KL}, \cite{St1}, \cite{DN}, \cite{DN2}, \cite{DKh} and \cite{DKh2},   
becomes useful for this purpose:  given an idempotent pure congruence $\rho $ on an inverse semigroup $S, $ there is a partial action $\tau $ of $S/\rho $ on the semillatice of the idempotents of $S$ such that $S$ embeds into the partial semidirect product 
determined by $  \tau $ (see    \cite[Theorem 3.4]{Khry1}). Recall that a  partial action of an inverse semigroup is   defined in  \cite{Khry1} by means of an arbitrary premorphism   into the symmetric inverse monoid, so that this notion is weaker than that one adopted in       \cite{LawMarSte} and
\cite{BuE}, where the premorphism is assumed to be order-preserving.   The image of $S$ in the partial semidirect  product is specified with the help of  the concept of a {\it fully strict partial action} of an inverse semigroup on a semilattice. 
It is also proved, using globalization, that O'Carroll's \cite[Theorem 4]{O'Carroll} is 
\cite[Theorem 3.4]{Khry1} with globalizable $\tau. $ It is pointed out by means of an example that $\tau $  is not   always globalizable.\\

Having elaborated the bases  of a cohomology theory of partial $G$-modules, one may try to extend the technique from 
\cite{DES2} to globalize partial $n$-cocycles. The main  result in \cite{DES2} asserts that an arbitrary  unital twisted partial action $\alpha $ of a group $G$ on a (unital) ring  ${\mathcal A},$ which  is a  product of indecomposable rings (blocks), admits an enveloping action, i.e. there exists a twisted global action $\beta $ of $G$ on a ring ${\mathcal B}$ such that ${\mathcal A}$ can be embedded into ${\mathcal B} $ 
as a  two-sided ideal, such that  $\alpha $ can be seen as  the restriction of $\beta $ to ${\mathcal A}$ and 
${\mathcal B}=\sum_{g\in G}{\beta }_g({\mathcal A})$. Moreover, if ${\mathcal B}$ has $1_{{\mathcal B}}$, then any two globalizations 
of $\alpha $  are equivalent in a natural sense.\\

If ${\mathcal A}$ is commutative, then the above mentioned results from \cite{DES2} mean  that given a unital 
$G$-module structure on ${\mathcal A}$, for any $2$-cocycle of $G$ with values in ${\mathcal A}$ there exists a (usual) 
$2$-cocycle $u$ of $G$ related to the global action on ${\mathcal B}$ such that $w$ is the restriction of $u$.  Moreover, if 
${\mathcal B}$ has $1_{{\mathcal B}}$, then any two globalizations of $w$ are cohomologous.\\

The proofs in \cite{DES2} are rather technical.  Nevertheless, thanks to some improvements  we managed   in \cite{DKhS} to extend to  arbitrary $n$-cocycles the results from \cite{DES2} in the commutative case. As in \cite{DES2} the globalization accours in two steps.  First we show that  given a unital partial $G$-module structure on a commutative ring 
${\mathcal A}$, a partial $n$-cocycle $w$ with values in ${\mathcal A}$ is glo\-ba\-li\-za\-ble if and only if a certain extendibility property holds for $w.$ The second, and more technical step, consists of establishing the extendibility property, assuming that  
 ${\mathcal A}$ is a product of blocks.\\

With respect to the  uniqueness of a globalization  we prove that  given a  glo\-ba\-li\-za\-ble unital partial action $\af$ of $G$ on a ring ${\mathcal A},$ such that ${\mathcal A}$ is a pro\-duct of blocks,  any two globalizations of   a partial $n$-cocycle $w$ related to $\af$ are cohomologous. More generally, arbitrary globalizations of 
cohomologous partial $n$-cocycles are also cohomologous. This gives an improvement even for the case $n=2$, as we do 
not require that $\af$ is unitally globalizable. This means that the global $n$-cocycles take values in the group of the 
invertible elements ${\mathcal U} ({\mathcal M} ({\mathcal B}))$ of the algebra   $ {\mathcal M} ({\mathcal B})$ of the multipliers of ${\mathcal B}.$ One  should notice that the global action of $G$ on  ${\mathcal B}$ naturally extends to    $ {\mathcal M} ({\mathcal B}).$
This allows us to  
establish an isomorphism between the partial cohomology group $H^n(G,{\mathcal A})$ and the global one 
$H^n(G,{\mathcal U} ({\mathcal M}( {\mathcal B})) )$.\\

Product of blocks were also considered in the more recent preprint \cite{CoMarcos}, in which the authors characterize partial actions of groups on a finite pro\-duct of indecomposable rings, dealing also with their enveloping actions.\\

An interesting recent advance on the globalization problem was obtained in  \cite{CFMarcos}, where the concept of a 
partial action of a group on a small weak $\kappa $-category was defined, as well as that of the corresponding partial skew category.
The latter  is a small weak  non-necessarily associative category, a  notion which is  also defined in   the same paper  \cite{CFMarcos}. 
The relation of the  partial skew category with  the concept of the  $\kappa $-algebra of a category and that of the category of a 
$\kappa $-algebra (with a given set of idempotents) is discussed, as well as the asocitivity of the    partial skew category.   
  Using the notion of an ideal in a category the authors define the concept of a restriction of a global action resulting in a partial action, as well as that of a globalization. Amongst other results, a criterion is given for the existence of a globalization for a partial action 
$\alpha $ of a group $G$ on a  small weak $\kappa $-category, under the assumption that $G$ acts globally on the objects via $\alpha .$ Moreover, if a globalization   exists, then  it is unique up to an equivalence.\\

Another latest   development  on the globalization problem was obtained in \cite{PiUz1}, were under appropriate conditions  a   criteria for the existence of a metrizable globalization for a given continuous partial action of a separable metrizable group $G$ on  a separable metrizable space $X$ was given. If  $G$ and $X$ are both Polish spaces, then  the globalization is   a Polish space too.  The existence of a universal globalization for  continuous partial actions of a countable discrete group  on Polish spaces is also  discussed. This topic was further studied by the same authors in \cite{PiUz2}, where amongst other results it was  proved  that the enveloping space  (i.e. the space under the enveloping action) of a partial action of a Polish group  on a Polish space  is a standard Borel space.\\

The globalization problem from the universal algebra point of view was investigated in \cite{KhryNov}, in which   a reflector of a partial action is constructed in the corresponding subcategory of global actions, and the question when this reflector is a globalization is considered. In particular,  the notion of a partial action of a group on a relational system is introduced and it is  shown that it admits a universal globalization which is a reflector. Partial group actions on partial algebras are also defined, whose domains are assumed to be   relative subalgebras, and for such  partial actions  a  necessary and sufficient condition for the existence of a globalization  is given.  For partial actions on total algebras  the above mentioned  reflector is constructed and shown to be the universal globalization.  For algebras with identities  the desired  reflector is also produced, but it may not be a globalization.  A characterization of globalizable, in the corresponding variety, partial actions is proved, which is applied to give  an example of a non-globalizable partial action  by taking the variety of semigroups.\\

The authors relate  the globalization problem in a variety of algebras  to embeddings of generalized    amalgams into an algebra from the variety. In the case of the variety of groups this is a  well-know and highly interesting topic, which was considered for other varieties too.   With any partial action $\theta $  of a group on an algebra from a variety 
$V$   a generalized amalgam $A$ of $V$-algebras is associated, such that $\theta $ is  globalizable if and only if $A$ is embeddable into a $V$-algebra.   The globalization problem is also considered in \cite{KhryNov} for partial actions on semigroups whose domains are ideals.\\

 In the majority of  cases the concept of the globalization considered by the authors is based on the restriction process described at the beginning of this survey, with domains given by formula (\ref{restriction}). Nevertheless, it may happen that a  partial action $\alpha = \{   \alpha _g : X_{g\m} \to X_g   \, ,   g\in G  \},$  extends to a
global action $\beta$  in such a way that  one has only the inclusions   $X_g \subseteq  Y \cap \beta _g (Y),$ rather than the equalities (\ref{restriction}). Moreover, it may be useful to extend a partial action $\alpha $ on a ring $\A$ to a better partial action $\alpha ^* $ on a
larger ring $Q$ without assuming that $\A$ is an ideal in $Q.$ If $\alpha ^* $ happens to be globalizable (in the above canonical sense), then the passage from $\alpha $ to the  globalization of  $\alpha ^* $ is definitely an  inetesting tool.
These ideas lead to weaker versions of the notion of a globalization.\\

Such a situation appeared first in M. Ferrero's paper \cite{F}, according to which a global action $\beta $ of a group $G$ on a ring $\B$  is called a {\it weak globalization} of a partial  action $\alpha  = \{   \alpha _g : D_{g\m} \to D_g   \, ,   g\in G  \}$ of $G$ on a ring $\A$ if  there is a monomorphism $\varphi : \A \to \B$ of rings such that  $$\beta _g  \circ \varphi | _ {D_{g\m}} = \varphi \circ \alpha _g \;\;\;  \text{for all} \;\;\; g\in G.$$ The main result in \cite{F} says that 
any   partial action $\alpha $ of a group $G$ on a semiprime  ring $\A$ possesses a weak globalization.\footnote{In fact, it is assumed  in \cite{F} that $\alpha $ is proper, i.e.  each $D_g$ is  non-zero, however, later in \cite{F2} M. Ferrero observed that this condition is unnecessary for the proof of this result.}  The proof is a  fruitful idea which was used in other articles in similar situations, the main step being an extension of $\alpha $ to a unital (and therefore globalizable)  partial action $\alpha ^*$ of $G$  on the Martindale ring ${\mathcal Q}$ of right quotients of $\A.$ Then the (canonical) globalization $(\B, \beta)$ of $\alpha ^*$ gives  a weak globalization of $\alpha .$ Notice that one does not need to assume that $\A$ has $1.$\\

 A subsequent article \cite{CF2} also deals with non-necessarily unital   semiprime rings, but the problem was considered from the point of view of the ca\-no\-ni\-cal globalization. A closure condition was imposed on the domains and   a property involving multipliers was used which appeared in \cite{DdRS} with respect to the globalization problem on $s$-unital rings. As we mentioned already, a partial actions  $\alpha = \{   \alpha _g : D_{g\m} \to D_g   \, ,   g\in G  \},$   of a group $G$ on a unital ring $\A$ is globalizable if and only of $\alpha $ is unital, i.e. each ideal $D_g$ is a unital ring. Now, in the case of a left $s$-unital $\A$ it is necessary, but not sufficient to assume that each $D_g$ is a left $s$-unital ring (i.e. $\alpha $ is left $s$-unital).     The criterion given in \cite{DdRS} says that $\alpha $ is globalizable if and only of $\alpha $  is left $s$-unital and for  each $g\in G$ and $ a \in \A$ there exists a multiplier $ \gamma _g(a) $ of $\A$ such that  

\begin{equation}\label{multiplier}
\A \gamma _g (a) \subseteq \A \;\;\text{ and } \; \;   x \gamma _g (a) = \alpha _g (\alpha \m _g (x) a ), \;\; \forall x \in  D_g.
\end{equation}
\noindent Now the main result in \cite{CF2} says that a partial action $\alpha $ of a group $G$ on a non-necessarily unital 
semiprime ring $\A,$ such that each ideal $D_g$ is closed,  is globalizable if and only if the above condition (\ref{multiplier}) is satisfied. Factoring the ring under the global action by the prime radical the authors come to a semiprime globalization, which is shown to be unique. A certain relation with the weak globalization from \cite{F} is discussed.\\ 

In a more recent article \cite{BeF} the above globalization result from \cite{CF2} was refined as follows.  Let 
  $\alpha = \{   \alpha _g : D_{g\m} \to D_g   \, ,   g\in G  \}$ be a partial action of a group $G$ on a non-necessarily unital semiprime ring $\A$ and    $(\B, \beta)$ be  the  weak globalization of $\alpha $ as above. Identifying $\A$ with its copy in $\B,$ write $W=\sum_{g\in G} \beta (\A),$ and let $\beta ' $ be the restriction of $\beta $ to $W.$ Observe  that $W$ is a ring and $\A$ is an ideal in $W$ if and only if $\A \beta _g (\A ) \subseteq \A,$ for all $g \in G.$ Producing multipliers which satisfy  (\ref{multiplier}) the authors prove that $(W, \beta ')$ is a globalization of  $\alpha ,$ provided that each $D_g$ is a direct summand of $\A.$\\  

It is essentially more complicated  to produce a weak globalization in the twisted case, which   was  considered in \cite{BeCoFeFlo} with more involved use of the quotient rings technique.  Let $\alpha $ be a twisted partial action of a group $G$ on a non-necessarily unital semiprime ring $\A.$ The authors in  \cite{BeCoFeFlo} extend first $\alpha $ (together with its mupliplier valued twisting)  to a unital twisted partial action $\alpha ^*$ on the left  Martindale ring of quotients $Q(\A)$ of 
$\A.$   Unfortunately, there are no known   tools to guarantee that $\alpha ^*$ is globalizable, unless $Q(\A)$ is a pro\-duct of indecomposable rings \cite{DES2}. So the authors go further and   extend  $\alpha ^*$ to a unital twisted partial action ${\alpha}^{**}$ of $G$ on the left maximal ring of quotients $Q_m(Q(\A))$ of $Q(\A).$ Actually, one may assume that  $Q_m(Q(\A)) = Q_m (\A).$  Now assuming in addition that $\A$ is a left Goldie ring, one has that  $Q_m (\A)$ coincides with the classical ring of quotients of $\A,$ which is semisimple by Goldie's Theorem. Consequently,   $Q_m (\A)$ is a direct product of in\-de\-com\-po\-sa\-b\-le rings and the already mentioned fact from \cite{DES2} on the globalization of twisted partial actions is applicable, resulting this way in a weak glo\-ba\-li\-za\-ti\-on for $\alpha. $\\

Ferrero's technique \cite{F}   was  used again in \cite{Avi1} to extend a unital (= glo\-ba\-li\-za\-ble) partial action $\alpha $ on an $\alpha $-semiprime (unital) ring $\A$ to the Martindale ring ${\mathcal Q} $ of $\alpha $-quotients of $\A.$ Denoting by 
$(\B, \beta)$ the (usual)
glo\-ba\-li\-za\-ti\-on of $(\A, \alpha ),$ and extending $\beta $ to an action $\beta ^*$ of $G$ on the  Martindale ring $Q  $ 
of $\beta  $-quotients of $\B,$ the author studies the relations between the involved actions and rings.  In particular, a criterion is given when  $\beta ^*$ is a globalization for $\alpha ^*.$\\

  The above mentioned paper \cite{DFP} extends to the partial action setting  the Galois Theory of commutative rings 
by   S. U.~Chase, D.K.~Harrison and  A.~Rosenberg \cite{CHR}, including several equivalent definitions of a partial Galois 
extension and  establishing a fundamental theorem on Galois correspondence. In  \cite{CHR} the authors also gave the 
exact sequence\\

$0\to H^1(G, \U(\A)) {\to} {\bf Pic}(\A^G){\to } {\bf Pic}(\A)^G {\to }H^2(G, \U (\A)) {\to} B(\A/\A ^\af){\to}$\\

\hspace*{3mm}$\to H^1(G,{\bf Pic}(\A )) {\to} H^3(G, \U(\A )),$\\

 \noindent which  generalizes the two most fundamental facts from Galois cohomology of fields, namely, the Hilbert's 
Theorem 90 and the isomorphism of the relative Brauer group $ B(\A /\A ^\af)$ with the second cohomology group 
$H^2(G, \A ^*).$ The latter isomorphism is obtained by associating to any $2$-cocycle from 
$Z^2(G, \A ^*)$ the corresponding crossed product $ \A  \ast G.$ The above sequence was  derived in \cite{CHR} from the Amitsur cohomology seven terms exact sequence by S. U. Chase and  A. Rosenberg \cite{CR},  specifying it  to the case of a Galois extension.  The proof in \cite{CR} used spectral sequences and was not constructive.  The first constructive proof was given by T. Kanzaki \cite{Kanzaki1968}, introducing and applying   generalized crossed products. Since then much attention have been payed to the  sequence and its parts establishing  
ge\-ne\-ra\-li\-za\-ti\-ons and analogues in various   contexts.\\

   Having developed the  partial Galois theory in \cite{DEP} on one hand, and the partial group cohomology in \cite{DKh} on the other, it was reasonable to establish the  analogue of the Chase-Harrison-Rosenberg  exact  sequence in the context of a 
partial Galois extension of commutative rings. The treatment in the partial action setting turned out to be more laborious, and 
some conceptual adjustments were needed to be made. The corresponding homomorphisms were constructed in \cite{DoPaPi1},
whereas the exactness of the sequence is proved in \cite{DoPaPiRo}. Amongst the new ingredients   we introduce the Picard monoids  ${\bf PicS}_{\A ^\af}(\A )$ and    ${\bf PicS} (\A ),$ the latter being  an inverse semigroup which is a disjoint union of the Picard groups
 of all direct summands of  $\A.$  Moreover,    a partial action $ \af^*$  of the Galois group $G$ on    ${\bf PicS} (\A )$ is used, 
as well as a partial action version of the generalized crossed products and two partial representations  of the form $G\to{\bf PicS}_{\A ^\af}(\A )$. As the final result we have obtained the following  seven-terms 
exact sequence\\

$0\to H^1(G,\af, \A ) {\to} {\bf Pic}(\A ^\af) {\to }{\bf PicS}(\A )^{\af^*}\cap {\bf Pic}(\A )  {\to }H^2(G,\af, \A ) {\to}$\\

\hspace*{2mm} $ \to B(\A /\A ^\af) {\to} H^1(G,\af^*,{\bf PicS}(\A )) {\to}  H^3 (G,\af, \A ).$\\

Other recent Galois theoretic results were produced  for partial group actions on rings in \cite{JiangSzeto2014}, \cite{JiangSzeto2015}, \cite{KuoSzeto}, \cite{KuoSzeto2} and \cite{KuoSzeto3}, and  for partial coactions on coalgebras in \cite{CasQua2}. The latest survey by A.\ Paques \cite{Paques2018} gives an overview of Galois theories, inlcuding those based on partial actions.\\

Partial representations are in the origin of the successful approach to the study of $C^*$-algebras generated by partial 
isometries via  partial actions.  By a partial isomtery in a $*$-algebra $\mathcal A$ 
(in particular, in a  $C^*$-algebra) one means an element $s\in \mathcal A$ with $s s^* s =s.$ By a projection in 
 $\mathcal A$ one understands a $*$-symmetric idempotent, i.e. an element $p\in   \mathcal A$ such that 
$p^* =p$ and $p^2=p.$ A {\it partial $*$-representation} $u: G\to  \mathcal A ,$ $g \mapsto u_g,$ $(g\in G),$  by 
definition is a partial representation such that $u_{g\m}= u_{g}^*,$ for all $g\in G.$ It follows from the definition that 
 each $u_g,$ $(g\in G)$ is partial isometry.\\

A product of partial isometries is not necessarily  a partial isometry, and an algebra generated by partial isometries, in 
general,  
may be rather wild,  unlikely to yield to any attempt at understanding its structure.   Nevertheless, algebras (abstract or $C^*$) generated by the range of a partial representation  have a chance to be endowed with the structure of a crossed product by a partial action, permitting one to understand  their  algebraic behavior.  One of the  first prominent examples of the use of this technique 
was established in the case of the Cuntz--Krieger algebras ${\mathcal O}_A$ defined as follows \cite{CK}:
given an $n\times n$ matrix $A = \{a_{ij}\}_{1\leq i,j\leq n}$ with
entries in $\{0,1\}$ one defines ${\mathcal O}_A$ as being the universal
$C^*$-algebra generated by partial isometries $S_1,\ldots,S_n$ subject
to the conditions:\smallskip

\begin{itemize}
\item[{CK$_1$)}] $\displaystyle\sum_{i=1}^n S_i S_i^* = 1$, and
  \smallskip 
\item[{CK$_2$)}] $\displaystyle S_i^* S_i = \sum_{j=1}^n a_{i,j} S_j S_j^*$.
\end{itemize}
  \smallskip

\noindent It was shown for them in 
\cite{E2} that there exists a partial representation
of the free group $\F_n$ sending the $i^{th}$ canonical  generator of  $\F_n$ to
$S_i$.  This idea was subsequently generalized in \cite{EL} to treat
the case of infinite matrices and was used to give the first
definition of  an analogue of Cuntz--Krieger algebras for
transition matrices on infinitely many states,  dropping the row-finite
condition used in  earlier investigations. The algebras defined in \cite{EL} are called the Exel-Laca algebras.\\

The key idea is as follows. Let  $u: G\to  \mathcal B $ be a partial representation into an algebra 
$  \mathcal B.$ If $\mathcal B $ is a $*$-algebra, then $u$ is assumed to be a partial $*$-representation. 
It is a well-known to the experts  fact that the  $ e_g = u_g u_{g\m},$ $( g\in G),$ form  a commutative set $E$ of idempotents,
which are projections in the $*$-case.    Let $  \mathcal A $ be the subalgebra of $\mathcal B $ generated by  $E.$ In the $C^*$-case  $  \mathcal A $ is the  $C^*$-subalgebra geberated by $E,$ i.e. the smallest   $C^*$-subalgebra of  $  \mathcal A $ containing $E.$ Let 
$D_g$ be the ideal in $  \mathcal A $ generated by $e_g$, i.e. $D_g = \mathcal A e_g, $ and let 
$\tau _g : D_{g\m} \to D_g$ be the map defined by 
\begin{equation}\label{pRep->pAc}
\tau _g  (a) = u_g a u_{g\m}, \;\;\; \forall g\in G.
\end{equation}
Then $\tau = \{ \tau _g : D_{g\m} \to D_g , g\in G \} $ is a partial action of $G$ on   $  \mathcal A $ 
\cite[Lemma 6.5]{DE} (in 
the $C^*$-case $\tau $ is a $C^*$-algebraic partial action, which can be seen by carrying over the purely algebraic proof). 
Furthermore, in the ring theoretic case by \cite[Proposition 6.8]{DE} there is a homomorphism from the crossed product (skew group ring) ${\mathcal A} \rtimes _{\tau} G $ to   $  \mathcal B, $ 
which  is an epimorphism in our case, as  $\mathcal B $ is generated by the elements $u_g, g\in G .$ In the $C^*$ case, 
thanks to \cite[Proposition 11.14]{E6}, there is a $*$-homomorphism  from the (full) $C^*$-algebraic crossed product 
${\mathcal A} \rtimes _{\tau}^{\rm full} G$ onto   
$\mathcal B .$ In some cases one is able to prove that the epimorphism obtained this way is, in fact, an isomorphism. Here we temporarily use the non-standard notation ${\mathcal A} \rtimes _{\tau}^{\rm full} G$ for the (full) $C^*$-algebraic crossed product
in order to make difference with the ring-theoretic crossed product ${\mathcal A} \rtimes _{\tau} G .$ One 
should also note that in the theory of  $C^*$-algebras there is also the so-called  reduced $C^*$-algebraic crossed product 
${\mathcal A} \rtimes _{\tau}^{\rm red} G$  by a partial action (see \cite{E6} for details).\\

A very recent use  of this technique in \cite{DE2} made it possible  to endow the Carlsen-Matsumoto $C^*$-algebra  
${\mathcal O}_X$ of an arbitrary   subshift $X$  with the structure  of a $C^*$-crossed product (in this case the full  and the 
reduced  $C^*$-algebraic crossed products coincide).  The approach based on partial actions and partial representations 
results in an alternative definition of ${\mathcal O}_X$, which is more convenient for our technique.\\

In order to describe briefly the  idea of this application, let  $\Lambda $ be  a finite alphabet and  $\Lambda ^{\N}$ be the set of all 
infinite words $ x_1x_2x_3\ldots $ in alphabet   $\Lambda $ (i.e. $x_i \in   \Lambda $). Taking the discrete topology on $\Lambda $ and the product topology on   $\Lambda ^{\N},$ one has that    $\Lambda ^{\N}$ is a  Hausdorff compact 
totally disconnected space. Then the (left) {\it shift map} $$S: \Lambda ^{\N} \to \Lambda ^{\N}$$ is defined by  
removing 
the first 
letter, and it  is easily seen to be con\-ti\-nu\-o\-us.   By a  (left) {\it subshift} $X$ one means  a closed $S$-invariant 
subset of     $\Lambda ^{\N}.$ Important examples of subshifts are obtained as follows. Given an arbitrary subset 
${\mathcal F } $
of finite words in alphabet $\Lambda$,  called the set of {\it forbidden words}, denote by
  $$
  X = X_{\mathcal F }
  $$
   the set of all infinite words $x$ such that  no member of $\mathcal F $ occurs in $x$ as an interval (contiguous 
block of letters). Then it is well known that $ X_{\mathcal F }$ is a subshift and any subshift is of this form. If $\mathcal F $ is 
finite, then 
 $ X_{\mathcal F }$ is called  {\it a subshift of  finite type}.
Subshifts are the objects of study of Symbolic Dynamics (see \cite{LindMarcus}), and it is interesting to bring new point of view 
relating algebras to them.\\

 It turns out that there is a natural partial action on $X,$ which is defined as follows.  For each  $a \in \Lambda$ denote by $X_a$ the set of all words in $X$ which begin with $a.$ Restricting the shift map $S$ we 
obtain the function  $X_a \to X,$ whose image will be denoted by $X_{a\m}.$ Evidently,   $X_{a\m}$ consists of all 
 $x\in X$ for which $ax \in X.$ Thus the restriction of $S$ results in the partial bijection $X_{a} \to X_{a\m},$
which will be denoted by $\theta_{a\m}.$ We also  write $\theta_a : X_{a\m} \to X_{a}$ for its inverse.\\

Let now $\F =\F(\Lambda )$ be the free group over $\Lambda ,$ i.e. the free group  whose set of canonical  generators is 
$\Lambda.$ Thus with each  generator 
$a \in \Lambda$ of $\F _n$ we may associate a partial bijection $\theta_a : X_{a\m} \to X_{a}$ of $X.$ Composing the
$\theta_a$'s  and their inverses one may associate to each $g\in \F _n$ a partial bijection 
$\theta _g : X_{g\m } \to X_g$ of $X,$ so that $\theta = \{ \theta _g, g\in G\}$ will be a partial action   of  $\F _n$ on $X$
(see \cite[Proposition 4.10]{E6}). We call $\theta$ the {\it standard partial action} associated to $X.$ In particular, any finite 
word $\gamma $ in alphabet $\Lambda$ is an element in   $\F _n$  and 
$$\theta _{\gamma } (x) = \gamma x , \;\; \theta_{\gamma \m }(\gamma x) = x \;\; \;\forall x\in X_{{\gamma}\m}.$$

Now taking the relative topology on $X$ one may wonder whether or not $\theta $ is a topological partial action. 
It follows by \cite[Proposition 2.5]{DE2} that $\theta$ is topological if and only if $X$ is of finite type. So we clearly have a difficulty with $\theta $ in the non-finite type case.  Subshifts of finite type
can be recoded and seen as Markov subshifts, i.e.  subshifts obtained from graphs (see \cite{LindMarcus}). Since 
the $C^*$-algebras related to Markov subshifts are the Cuntz-Krieger algebras, which are well-understood from the point of view of partial actions, we are mainly interested in the case of subshifts of non-finite type, so the above mentioned problem 
with $\theta $ should be resolved somehow.\\

Despite the bad topological behavior of the standard partial action,  we may use it to define a partial $*$-representation 
of $\F $ by bounded operators, which will be crucial for our approach.  For each  $g$ in $\F $, 
denote by $u _g$  the unique bounded linear operator
  $$
  u _g: \ell ^2(X) \to \ell ^2(X),
  $$
  such that for each  $x$ in $X$,

$$   u _g(\delta _x) =\left\{
\begin{array}{ll} \delta _{\theta _g(x)}  , & \text{if } x\in X_{g\m },\\
 0,  &  \text{otherwise},
\end{array} \right. $$ where $\theta $ is the standard partial action associated to $X$.
  Then the map $g\mapsto u _g$ is a  partial $*$-representation of  $\F $ on $\ell ^2(X)$.\\

Denote by ${\mathcal M}_X$ the closed $\ast $-algebra of bounded operators on $\ell ^2(X)$ ge\-ne\-ra\-t\-ed by 
$\{u_g \, : \, g\in \F \}.$ It follows by \cite[Proposition 6.1]{DE2} that ${\mathcal M}_X$ coincides with the $C^*$-algebra 
defined by K. Matsumoto in \cite[Lemma 4.1]{MatsuAutomorph} (see also \cite{CarlMatsu2004}), and we call   
${\mathcal M}_X$ the {\it Matsumoto algebra} associated to $X.$\\

Denoting by  $\lambda $   the left regular representation of the free group $\F =\F (\Lambda )$ on ${\ell} ^2(\F )$ and   
by  $\{\delta _g\}_{g\in \F }$  the canonical
orthonormal basis of ${\ell} ^2(\F ),$ one has that
  $$
  \lambda _g(\delta _h) = \delta _{gh},
  \forall g,h\in \F .
  $$

Using  the partial representation $u $  we  define a new partial representation
$\tilde u $ of $\F $ on ${\ell} ^2(X)\otimes {\ell} ^2(\F )$, by tensoring $u $ with $\lambda $, namely
  $$
  \tilde{u} _g = u _g\otimes \lambda _g, \forall g\in \F .
  $$

\noindent Then we define the {\it Carlsen-Matsumoto $C^*$-algebra} ${\mathcal O}_X$ associated to a given subshift $X$ as the closed
$\ast$-algebra of
operators on $\ell ^2(X)\otimes \ell ^2(\F )$ generated by the set $\tilde{u} (\F ) = \{ \tilde{u}_g \, g \in \F\}.$ It follows  from \cite[Theorem 10.2]{DE2} that the above defined  $C^*$-algebra ${\mathcal O}_X$ is isomorphic the    $C^*$-algebra defined by 
T.~M.~Carl\-sen in \cite[Definition 5.1]{CarlsenCuntzPim}.\\

Having at hand the partial representation $\tilde{u}: \F \to  {\mathcal O}_X$ we use the above mentioned technique  (\ref{pRep->pAc}) to 
produce a $C^*$-algebraic partial action $\tau $ of $\F $ on the commutative  $C^*$-subalgebra ${\mathcal D}_X$ of  ${\mathcal O}_X$ generated by the projections $ \tilde{u} _g\tilde{u} _g^*,$ $g\in \F.$ One readily identifies  
${\mathcal D}_X$ with the $C^*$-subalgebra  of ${\mathcal M}_X$ generated  by the projections $ e_g =\tilde{u} _g\tilde{u} _g^*,$ $g\in \F,$ as  $ \tilde{u} _g\tilde{u} _g^* = e_g \otimes 1.$ Then it is straightforward  to 
  shows that there exists a  $*$-homorphism
 ${\mathcal D}_X \rtimes _{\tau}^{\rm full} \F \to {\mathcal O}_X$ \cite[Proposition 9.2]{DE2}, which is surjective, as  
${\mathcal O}_X$ is generated by the elements $\tilde{u}_g,$ $g\in \F.$ In fact, further considerations and known results
lead to the following isomorphisms of $C^*$-algebras:

$$ {\mathcal D}_X \rtimes _{\tau}^{\rm full} \F \cong {\mathcal O}_X  \cong {\mathcal D}_X \rtimes _{\tau}^{\rm red} \F.$$

\vspace*{2mm}

By Gelfand's Theorem,   the $C^*$-algebra ${\mathcal D}_X$ is isomorphic to $C (\Omega _X),$ where $\Omega _X$
stands for the spectrum of ${\mathcal D}_X.$ It is known by the $C^*$-theory of partial actions (see [Corollary 11.6]\cite{E6}) that there is a  topological partial action $\vartheta $ of $\F $ on  $\Omega _X$ which corresponds to $\tau ,$ and which we 
call {\it the spectral partial action}. Using $C^*$-theoretic notation the above isomorphisms may be rewritten in the form

$$ C(\Omega _X) \rtimes _{\vartheta}^{\rm full} \F \cong {\mathcal O}_X  \cong 
C({\Omega }_X ) \rtimes _{\vartheta}^{\rm red} \F.$$

\vspace*{2mm}

Due to the above mentioned difficulty with  the standard partial action, the main topological partial dynamical system 
related to an arbitrary subshift is $(\Omega, \vartheta),$ and for this reason a special attention should be payed to the topological space $\Omega.$ In \cite{DE2}  some  properties of the elements of   $\Omega $ are given,  so that 
if we look at them as subsets of the Cayley graph of $\F ,$ then they have  the  aspect of a river basin.  It does not seem 
to be 
 possible to  find a complete set of properties characterizing $\Omega $, but we  know that it contains a dense
copy of $X$ (not necessarily with the same topology), permitting to deal with  the other,   more elusive elements. We are able 
to  give necessary and sufficient conditions for such relevant properties of $\vartheta $ as minimality and topological 
freeness, in a rather ``graphical'' language,  similar to those known for the case of Markov subshifts 
\cite[Theorems 12.6,  13.19]{E2}.  Then we use these
conditions to give a criterion for the simplicity of $\mathcal O _X$ in terms of $X$ \cite[Theorem 14.5]{DE}.\\ 

The above shows the importance of  a partial representation $u: G\to \A$  when the algebra $\A $ (abstract or $C^*$)  is generated  by the range $U$ of $u.$ If  $\A $ is a $*$-algebra then in this case we say that $\A $ is generated by a {\it tame} 
set of partial isometries. Equivalently, each element of the multiplicative semigroup $\langle U \cup U^* \rangle$ generated 
by $ U \cup U^*$ is a partial isometry   (see   \cite[Definition 12.9 and Proposition 12.13]{E6}).
Even if a $*$-algebra is generated by a wild (non-tame) set of partial isometries, one may impose additional relations on the generators in order to force them to behave well. This is concretely done in the recent paper
\cite{AraEKa} for the case of the Leavitt $C^*$-algebra $L_{m,n}.$ The latter is the topological analogue of the algebra 
$L(m,n)$ considered  by W. G. Leavitt in \cite{Leavitt62}.  The case $L(1,n)$ is a part of the well developed theory of the 
Leavitt-path algebras, which are ring-theoretic analogues of the graph $C^*$-algebras.  The generating set $U$ of 
$L_{m,n}$ is a wild set of partial isometries, and a way to turn around of this difficulty is to consider the  $C^*$-algebra 
 ${\mathcal O}_{m,n}$ which is the quotient of $L_{m,n}$ by the closed  ideal generated by the elements $x- x x^* x,$ 
where $x$ runs over  the semigroup $\langle U \cup U^* \rangle .$  This guarantees the tameness of the 
generating set of  ${\mathcal O}_{m,n},$ allowing one to apply the above mentioned technique based on partial representations and partial actions. As a consequence, 
${\mathcal O}_{m,n}$ is shown to have a full $C^*$-algebraic crossed product structure by a partial action (see \cite[(2.5)]{AraEKa}). It is 
in\-te\-res\-ting to notice that, unlike  the case of algebras related to subshifts, 
${\mathcal O}_{m,n}$ is not isomorphic to the corresponding reduced  $C^*$-algebraic crossed product 
\cite[Theorem 7.2]{AraEKa}.\\   

Analogously to the case of (usual) representations, there is an algebra responsible for the partial representations. In the 
case of the $*$-representations of a group by bounded operators,  the {\it partial group $C^*$-algebra} $C^*_p(G)$ defined 
in  \cite{E1} plays this role, whereas its ring theoretic version $\kpar \, G$  is the semigroup algebra $\kappa  {\mathcal S}(G),$  considered    in \cite{DEP}. According to a structural result from \cite{DEP},  if $G$ is a finite group  and $\kappa $ is a commutative ring (which is assumed to be associative and unital), the 
the partial group algebra  
 decomposes as follows:
	\begin{equation*}
	\kpar \, G \cong \bigoplus_{
		\begin{tiny}
			\begin{array}{c}
				H\leq G \\
				1\leq m \leq \left[ G:H \right]
			\end{array}
		\end{tiny}}
	\frac{b_m(H)}{m}\, \mathbb{M}_m\left(\kappa H\right),
	\end{equation*}
where $\kappa H$ stands for the group algebra of $H$ over $\kappa $  and $b_m(H)$ denotes the number of subsets $A\subseteq G$, such that $|A|= m |H|,$ $1\in A,$ $H=\left\{g\in G \tq gA=A\right\}.$ A recursive formula for the coefficients $\frac{b_m(H)}{m}$
given in \cite[(14)]{DEP} was corrected in \cite[(2)]{DP}. In \cite[Theorem 2.4 ]{Choi1} K. Choi gave an interesting  formula for the  coefficients  $b_m(H)$ using the M\"obius function, but one needs to be careful with Choi's notation in  \cite{Choi1} (see \cite[Remark 3.2]{DS2}).\\

The proof of the above structural result used the finite groupoid  $\Gamma=\Gamma(G)$ associated to a finite group $G,$
the elements of   $\Gamma $  being the pairs $(A,g)$, where $g\in G$ and $A$ is a  subset 
of $G$ containing the identity $1=1_G \in G$ and the element $g^{-1}$. 
The product of pairs  $(A,g)\cdot(B,h)$
in $\Gamma$ is defined  for pairs for which   $A=hB$, in which
case:
$$(hB,g)\cdot(B,h)=(B,gh).$$ Then one considers   the groupoid algebra 
 $\kappa \Gamma(G),$ which is unital  with $$1_{\kappa  \Gamma(G)} = \sum_{A\ni 1_G}(A,1_G),$$ together with the 
map $\lpar:G\longmapsto \kappa \Gamma(G)$:
\begin{equation}\label{parlambda}
\lpar(g)=\sum_{A\ni g^{-1}}(A,g),
\end{equation} which turns out to be a partial representation.\\

It is readily seen that any partial representation $u : G \to \A $ into a unital algebra $\A$ uniquely extends to a 
homomorphism of algebras  $\kpar \, G \to \A$ given by $[g] \mapsto u(G),$ $(g\in G).$  Theorem 2.6 in \cite{DEP} says 
that a similar universal property holds for $\kappa  \Gamma (G),$ from which we derive that   $\kpar \, G \cong \kappa  \Gamma (G),$
provided that $G$ is finite. This identification with the  groupoid algebra  then easily leads to the above structural fact.\\

After the publication of \cite{DEP} it was mentioned in private communications  independently by several 
researchers\footnote{Eric Jespers,  Stanley Orlando Juriaans, Boris Novikov.}  familiar with   inverse semigroups  that  the above structural result for partial group algebras can be obtained  using the theory of semigroup algebras. We spelled out  this in detail in   \cite[Remark 2.5]{DdLP}, without  claiming  novelty neither for the idea nor 
for the proof.\\

It was observed by R.~Exel in \cite{E1} that the complex partial group algebras  $\Cpar \, Z_4  $ and $\Cpar \, [ Z_2\times Z_2]$ of the cyclic group $ Z_4 $ of order $4$ and the Klein-four group, respectively,  are not isomorphic, whereas their usual complex group algebras  are isomorphic to $\C ^4.$ This suggested to consider the isomorphism problem for the partial group algebras:

\begin{center}

 Which properties of $G$ are determined by  $\kpar \, G$?\\

\end{center}

 \noindent In particular, does the isomorphism of algebras
  $ \kpar \, G_1  \cong \kpar \, G_2$ imply $G_1 \cong G_2?$\\

 A negative answer for the latter question was given 
already in
\cite{DEP}, producing two non-isomorphic finite non-commutative groups of order $605$  with isomorphic complex partial group algebras.  However, it was proved in \cite[Theorem 4.4]{DEP}  that if $G_1$ and $G_2$ are finite abelian groups and 
$\kappa $ is an integral domain with ${\rm char} \kappa $ not dividing $| G_1 | =|G_2 |,$ then $\kpar \, G_1 \cong \kpar \, G_2$ exactly when $G_1 \cong G_2.$ An analogous result in the mo\-du\-lar case (i.e. when ${\rm char} \,\kappa  \;$ divides $ \;|G_1|$)  was established in \cite[Theorem 3.7]{DP}. The above mentioned counter-example from \cite{DEP}   for the isomorphism question was improved in \cite{DS} by giving 
counter-examples over $\C $ for orders $|G|=2^5$ and $|G|=3^5,$ and  pointing out that there are no counter-examples for $|G| <2^5$ and $|G|=p^n < 3^5, p\neq 2.$\\

Notice that in the above mentioned positive results  on the isomorphism question from \cite{DEP} and \cite{DP} it is assumed 
that both given groups $G_1$ and $G_2$ are abelian. So one may wonder whether or not the partial group algebra 
$\kpar \, G$ determines the commutativity of $G.$\\

 More specifically, suppose that $G_1$ and $G_2$ are finite groups and
$G_1$ is abelian. Is true that  
$$\Cpar \, G_1 \cong  \Cpar \, G_2 \;\; \Longrightarrow \;\; G_2\; \text{is abelian?}$$

In \cite{DS}  a list of integers was obtained which form a complete list of invariants of $\kpar \, G,$ where $G$ is a finite $p$-group and $\kappa $ is an algebraically closed field with ${\rm char} \, \kappa  \neq p.$ The invariants were expressed  as sums in terms of subgroups of $G$ and numbers of their irreducible $\kappa $-representations.   Then
 they  were applied to the isomorphism problem, and  it was shown, in particular, that $\kpar \,G$ determines the commutativity of a finite $p$-group $G$ where $p$ is an odd prime.\\

The list of the invariants includes, of course, the order of the group, and this does not depend on the fact of $G$ being a $p$-group.  Indeed, for an arbitrary finite group $G$ the isomorphism   $\kpar \,G \cong \kappa  \Gamma (G)$ implies that 
the dimension of    $\kpar \,G $  is equal to 
\begin{equation*}
{\rm dim}(\kappa \Gamma(G))=\sum_{k=0}^{n-1}(k+1)\binom{n-1}{k}=2^{n-2}(n+1), 
\end{equation*}  where  $\vert G\vert=n.$ The right hand side of the above formula  is a strictly 
increasing function on $n,$ so that   if $G_1$ and $G_2$ are finite groups
such that  $\kpar \, G_1 \cong \kpar \, G_2$, then $\vert G_1\vert=\vert G_2\vert$.\\

 It turns out to be much more  complicated to obtain a full list of  invariants for the case when $G$ is  an arbitrary finite group. Ne\-ver\-the\-less, in \cite{DS2}    a series of  natural invariants were given, which are useful for the isomorphism question, provided that the group order is odd.   We point out the following two of them:\\
   
\noindent 1) The number of Sylow $p$-subgroups of  $G$ (with odd $|G|$);\\

 \noindent 2) $[H:H']$ where $H$ is a Sylow $p$-subgroup of $G.$\\

 \noindent   In particular, it follows that   if  $|G|$ is odd, then  $\Cpar \, G$ determines the commutativity of the group 
$G$ \cite[Corollary 5.3]{DS2}.\\

The above mentioned commutativity problem is much more difficult  for finite groups of even order. In particular, the following question remains open.\\

\noindent {\bf Open problem:} Suppose that $G_1$ and $G_2$ are finite $2$-groups such that $G_1$ is abelian and 
$\Cpar \, G_1 \cong  \Cpar \, G_2.$ Is it true that $G_2$ is also abelian?\\

We were not able to find a counter-example for the above problem using GAP.\\

 The invariants in \cite{DS} were obtained using a  formula for  $b_m(H),$ which involves counting   in the subgroup lattice of  a finite $p$-group $G.$ The formula  was generalized in \cite[Theorem 3.1]{DS2} and used to give an alternative proof for the above mentioned Choi's formula  \cite[Theorem 2.4 ]{Choi1}. It was also proved 
 that if $G_1$ and $G_2$ are finite groups such that there exists an isomorphism between the subgroup lattices of $G_1$ and $G_2$ which preserves subgroup rings, then $\kpar \, G_1\cong\kpar \, G_2 ,$ where $\kappa $ is a commutative ring 
\cite[Theorem 4.2]{DS2}.  This improves an earlier  result from  \cite{DP} in which it was additionally assumed that the lattice isomorphism preserves subgroups conjugation.  Several examples were also given in \cite{DS2}, including a counter-example 
for the  isomorphism problem of partial group algebras over $\Q $ with $|G_1| =|G_2|= 243= 3^5$, as well as the least counter-example over $\C $  in the case of  odd group order   ($|G_1| =|G_2|= 189 = 3 \cdot 7 \cdot 9$).\\

The decomposition of  partial group algebras into a product of matrix algebras over subgroup rings  gives information about 
the structure of the  matrix partial representations of the finite group $G$.     If $G$ is infinite, no such structural result is known, in particular  $\kpar \, G $ is not isomorphic to $\kappa  \Gamma (G),$ as  $\kappa  \Gamma (G)$ is not unital. Nevertheless, it is possible to relate the finite dimensional representations of $\kpar \, G $ with those of  $\kappa  \Gamma (G).$ To this end, for each 
connected component $\Delta $ of $\Gamma (G),$ whose set $ \V_{\D}$ of vertices is finite,  one uses an analogue 
$\lambda _{\Delta} : G \to \kappa   \Delta, $ of the   partial representation (\ref{parlambda}).  The map $\lambda _{\Delta} $ is defined by
\begin{equation*}
\lambda _{\Delta} (g)=\sum_{\substack{A\in \V_{\D}\\ A \ni g^{-1}}}(A,g),
\end{equation*} 
\noindent and it is shown to be a partial representation too. Identifying $K \Delta $ with a matrix algebra of the form 
$M_n (\kappa H)$ for some subgroup $H$ of $G, $ the map $\lambda_{\Delta}$ becomes a partial representation of the form 
$G \to M_n(\kappa H),$ called {\it elementary}. Then the main result Theorem 2.2 of \cite{DZh} says that composing the elementary partial representations of $G$ with the irreducible (usual)  representations of  $\kappa  \Delta $ (for all such $\Delta $), we obtain all irreducible partial representations of $G.$  A similar fact  holds for the indecomposable partial representations, by replacing above  the term ``irreducible'' by  ``indecomposable''. The   proof is rather technical, and one of the steps asserts that 
the  subalgebra of $\kappa  \Delta $ generated by the elements $\lambda _{\Delta} (g),$ $g\in G,$ coincides with $\kappa  \Delta .$ The argument given in \cite{DZh} needed a correction which was done in \cite[Proposition 2.2]{DdLP}.\\

The above mentioned domains $X$  of the partial projective representations of $G$ where further studied in the recent article \cite{DdLP}, mentioned above.
It was shown already in \cite{DN} that they are exactly the  $\mathcal{T}$-invariant subsets  of $G\times G,$  where $\mathcal{T}$ is a semigroup of order $25$ acting on  $G\times G.$ Surprisingly, the structure of  $\mathcal{T}$ does not depend on $G$ and it is a disjoint union of the symmetric group $S_3$ and an ideal which is a completely $0$-simple semigroup.  The domains $X$    form a lattice  $C(G)$ with respect to the set-theoretic inclusion, intersection and union, and  we know from \cite[Theorem 4]{DN2} that in the case of of a finite group $G$ they  are exactly the finite unions of {\it elementary domains},  i.e.  domains of  elementary partial representations.\\ 

The main structural pieces in   $C(G)$ are the atoms (i.e. minimal domains) and the $\mathcal{T}$-orbits which are not atoms. The latter are called  {\it blocks}.  Corollary  3.18 from \cite{DdLP} establishes the uniqueness of a minimal decomposition of a domain into a union of atoms and blocks.  The domains are strongly related to the ideals of the semigroup $S_3(G)$ \cite[Proposition 5.3]{DoNoPi},
so that the above fact can be translated into the uniqueness of a decomposition of any non-zero ideal in $S_3(G)$  into an intersection  of meet-prime ideals. The main decomposition result in \cite{DdLP}  is Theorem 4.3 which   in the case of a finite $G$ gives a decomposition of a non-minimal elementary domain into blocks and asserts its uniqueness under the assumption of  minimality. On the other hand,    Co\-rol\-lary 4.8 in \cite{DdLP} states the corresponding fact for the ideals in    
$S_3(G).$\\ 

The concept  of a partial group representation was extended to  that of a partial representations of a Hopf algebra in  \cite{AB2} and  \cite{ABV}. In fact, the first version of a definition of a partial representation of a Hopf algebra was given in  the article \cite{AB2}, where  the globalization problem for partial Hopf actions on algebras was also studied. Influenced by \cite{CaenJan} the latter was considered with respect to restrictions of global actions on right ideals (contrary to the above mentioned idea of restrictions on two-sided ideals). This generality resulted in an economical but an asymmetric concept of a partial representation \cite[Definition 6]{AB2}.  In many important results on partial Hopf actions one assumes that the partial actions are  symmetric 
\cite[(PA4)]{ABV},  a condition which is always satisfied by partial group actions. Using symmetric partial Hopf actions one 
comes to a new concept of a partial representation of a Hopf algebra \cite[Definition 3.1]{ABV} which can be obtained from
the initial de\-fi\-ni\-ti\-on given in  \cite{AB2} by adding three more axioms.\\

An interaction between partial actions and partial representations, similar to that in \cite{DE} was also established in the Hopf case \cite{ABV}. Furthermore, in \cite{ABV} for any Hopf algebra $H$ a Hopf  analogue $H_{par}$ of     $\kpar \, G$ was also introduced and studied. 
In particular, the map $H \ni h \to [h]\in H_{par},$ where  the $[h]$'s are the defining generators of  $H_{par}$ (see \cite[Definition 4.1]{ABV}), is a partial representation, which factors any other partial representation of the Hopf algebra $H$ \cite[Theorem 4.2]{ABV}, similarly to the group case.   The above mentioned technique  (\ref{pRep->pAc}) to produce partial crossed products  from partial group representations permits one to endow  $\kpar \, G$ with a crossed pro\-duct structure ${\mathcal A} \rtimes _{\tau} G,$ where $\tau $ is a partial action of $G$ on the subalgebra $\A \subseteq \kpar \, G,$ generated by the idempotents  $[g][g\m],$ $(g\in G)$ 
\cite[Theorem 6.9]{DE}. Analogously,  $H_{par}$ is a partial smash product of the form $\underline{\A \# H}$ \cite[Theorem 4.8]{ABV}, where $\A $ is the subalgebra of $H_{par}$ generated by the elements $\sum_{(h)} [h_1][S(h_2)],$ $(h \in H).$ A universal characterization of  $\A $ is given in \cite[Theorem 4.12]{ABV} under the assumption that the antipode of $H$ is invertible.\\

The Hopf case turns out to be more complex, enjoying  features which were invisible in the group case.  In fact,  $H_{par}$
has the structure of a Hopf algebroid \cite[Theorem 4.10]{ABV}, resulting in a functor from the category of Hopf algebras to that of Hopf algebroids \cite[Proposition 4.11]{ABV}. Moreover, we know that if $G$ is a finite group then    $\kpar \, G$ is finite dimensional. In the Hopf case, however, Example 4.14 shows that   $H_{par}$ may be infinite dimensional, even if $H$ is finite dimensional (this is the case of the Sweedler Hopf algebra $H_4$).\\

The obtained facts lead to an isomorphism of the category of      $H_{par}$-modules and the category of the so-called {\it partial $H$-modules} (and naturally defined morphisms). Here the reader should notice a   difference between the use of the term ``partial module'' in \cite{DKh}, \cite{DKh2}, \cite{DKh3} and in \cite{ABV}.  As it was explained above, in the  partial cohomology theory of groups by a partial $G$-module we understand a commutative monoid $A$ with  a partial action of 
$G$ on $A,$ whereas  in \cite[Definition 5.1]{ABV},  for any Hopf algebra $H$ over a field $\kappa ,$ by a partial $H$-module  
the authors mean a $\kappa $-vector space $M$  together with a partial representation of the form $H\to {\rm End}_{\kappa}  (M).$ The latter extends to the Hopf case the notion of a $G$-space given in \cite{DZh}, and can be reformulated in terms of certain partial actions of $G$ on $M$ \cite[Example 5.4]{Ba}. Another categorical isomorphism is established between the category of symmetric partial $H$-module algebras and the algebra objects  in the category of partial $H$-modules \cite[Theorem 5.6]{ABV}.\\

The authors in \cite{ABV} also obtained several fact on the so-called partial $G$-gradings. It is well-known that $G$-gradings on algebras are exactly the coaction of the  group Hopf algebra $\kappa G $ on algebras. Thus one may define a {\it partial grading} of $G$ on an algebra $\A $ as a partial coaction of   $\kappa G $ on $\A.$ In particular, partial representations are applied to obtain a complete characterization of the partial 
$\Z _2$-gradings on algebras  \cite[Theorem 6.8]{ABV}, which, in other words, gives a transparent description of 
all partial superalgebras as a product of an arbitrary (usual) superalgebra and its ideal. Even earlier in \cite{AB1}, an example was given of a partial  action of the dual group algebra  $(\kappa G)^{\ast}$ on an ideal $I$ of the group algebra $\kappa G,$ which, of course,  can be considered as a partial $G$-grading on $I.$\\

More generally,  in \cite{AAB}  the authors introduce and study the concept of a partial action of a Hopf algebra $H$ on a $\kappa $-category, also called a {\it partial $H$-module category}, and consider  
group gradings on $\kappa $-categories  as a special case of  Hopf module categories.  Partial $H$-module categories generalize both partial Hopf actions on algebras and  $H$-module categories. The latter were introduced  in \cite{CibSol} in order to relate smash extensions to Galois coverings of $\kappa $-categories. The authors of 
\cite{AAB} prove that any  partial $H$-module category ${\mathcal C}$ possesses a unique, up to isomorphisms,   minimal globalizatrion ${\mathcal D}$, define the  notion of the partial smash product category $\underline{{\mathcal C}\# H}$ and prove that  $\underline{{\mathcal C}\# H}$ and  ${{\mathcal D}\# H}$ are  Morita equivalent (in the sense of  \cite{CibSol}). These facts generalize the analogous results obtained earlier in \cite{AB2} for Hopf actions on algebras. Furthermore, a structural result is established for the partial gradings of $\kappa $-linear categories by finite groups whose order is not divisible by ${\rm char}\, \kappa,$ and several examples are elaborated \cite{AAB}.\\

The symmetry between the algebra and coalgebra structures of a Hopf algebra, as well as  the concept of the finite dual of a Hopf algebra,    provides  a natural  ambient  for  
duality phenomena.  The duality between actions and coactions of locally compact groups on von Neumann algebras \cite{NakaTake} inspired 
M. Cohen and S. Montgomery   \cite{CohMon} to study duality between group actions and group coactions on algebras, which was later generalized by  R.J. Blattner and  S. Montgomery \cite{BlatMon} for the Hopf algebra framework. In  \cite{Lomp}  the Cohen-Montgomery   
theorem was extended for partial actions of groups.  Further generalizations were made in \cite{AB3} for partial Hopf 
actions.\\

  A comprehensive treatment of duality between partial  Hopf actions and  partial Hopf coactions was recently given in   \cite{BatiVerc1}. In particular,  under certain  (co)commutativity conditions
a connection between partial (co)actions and Hopf algebroids was established in    \cite{BatiVerc1},  one of the  important steps being 
to endow the partial smash product $\underline{\A \# H}$ with a structure of an  $\A$-Hopf algebroid, provided that $H$ is cocommutative and $\A$ is commutative. Furthermore,   for  a commutative Hopf algebra $H$ and a right partial $H$-coaction on a commutative algebra $\A$ the authors  construct a commutative $\A$-Hopf algebroid 
$ \A \underline{\otimes} H$, the so-called  partial split Hopf algebroid. In the opposite direction, from certain commutative $\A$-Hopf algebroids the authors come to partial Hopf coactions. Moreover,  a so-called skew-duality is established between   $\underline{\A \# H}$  and  
$\A \underline{\otimes} H,$  under certain conditions. In addition, considering the concept of a  partial Hopf action and that of a    partial Hopf coaction on coalgebras,  an interaction between the corresponding  $C$-rings  and  cosmash coproducts is also explored.  Several other valuable facts are also proved and interesting examples are given, one of them based on algebraic groups. Without giving formal definitions,  partial actions of affine algebraic groups were used earlier in \cite{ABDP} to produce concrete  partial Hopf coactions and then, dualizing, a partial Hopf action.   
Now the authors formalize the definition and show, as   expected,  that there is a one-to-one correspondence between the partial actions of an affine algebraic group $G$ on an affine variety $M$ and the partial Hopf coactions of the  coordinate Hopf algebra of $G$ on the coordinate ring of $M$ \cite{BatiVerc1}, which generalizes a similar classical fact known in  the global case.\\ 

   Hopf algebroids appear also  with respect to  partial Hopf cohomology,   which is the matter of the hot off the press preprint  \cite{BaMoTe}. As it was mentioned above, partial $G$-modules in   \cite{DKh} mean unital partial actions of a group $G$ on a commutative monoid  $A.$ Assuming that $A$ is a (unital) algebra, it is possible to extend the notion of a partial cohomology group  from   \cite{DKh} to the context of a partial action of  a co-commutative Hopf algebra $H$  on a commutative algebra $A$ \cite{BaMoTe}, giving  a natural generalization of Sweedler's cohomology  \cite{Sweedler} to the partial action setting. In    \cite{DKh} the monoid $A$ is replaced by an appropriate inverse submonoid $\tilde{A},$ which does not change the cohomology. In the Hopf case   $\tilde{A}$ is a certain quotient algebra of the symmetric algebra generated by the images of all reduced partial $n$-cochains.  The algebra $\tilde{A}$  possesses a partial action of $H,$ so that, as in the group case,
the isomorphism $H^n_{par} (H, A) \cong  H^n_{par} (H, \tilde{A})$ of the corresponding co\-ho\-mo\-lo\-gy groups holds. Furthermore, $ \tilde{A}$ enjoys a structure of a commutative and co-commutative Hopf algebra 
\cite[Theorem 4.5]{BaMoTe}.  In addition, by a result from \cite{ABDP} one naturally concludes that the partial crossed products 
$\underline{A \# _{\omega} H}$ (with commutative $A$ and co-commutative $H$) are in a bijective correspondence with the cohomology classes  
$[\omega] \in  H^2_{par} (H, A).$\\

 A significant new feature in \cite{BaMoTe} is the fact that the partial crossed product $\underline{\tilde{A} \# _{\omega} H}$ is a Hopf algebroid over the base algebra  
$E(A) = \langle h \cdot 1_A \, | \, h\in H \rangle $ \cite[Theorem 5.10]{BaMoTe}. The latter suggests that  it should be possible to look at the  partial Hopf cohomology from the point of view of the cohomology of Hopf algebroids 
(see \cite{BohmStefan}). Another important Hopf algebroid comes into the picture with respect to the notion of a partially Cleft extension from \cite{ABDP}, which, in the present setting,  is shown to be  a cleft extension by a Hopf algebroid in the sense of \cite{BohmBrzezinski}. More precisely, given  a co-commutative Hopf algebra $H$ acting partially on a commutative and co-commutative Hopf algebra $\tilde{A},$ and  a partial 2-cocycle $\omega \in Z^2_{par} (H, \tilde{A}),$  the partial smash pro\-duct  $\underline{E( \tilde{A} ) \#  H}$ is a Hopf algebroid over $E(  \tilde{A} )$ by \cite[Theorem 3.5]{BatiVerc1}. Then \cite[Theorem 6.6]{BaMoTe} states that  the partial crossed product 
 $ \underline{A \# _{\omega} H}$   is a right    $\underline{E( \tilde{A} ) \#  H}$-module algebra whose subalgebra of the co-invariants is $  \tilde {A} ,$ and the 
extension 
$\tilde{A} \subseteq \underline{A \# _{\omega} H}$ is a cleft extension by the Hopf algebroid  $\underline{E( \tilde{A} ) \#  H}.$ Recalling from \cite{ABDP} that any partially cleft $H$-extension over the co-invariants can be seen as a crossed product extension, we obtain the desired fact. The article is enriched with examples and suggestions for further research.\\

As we mentioned already, the term ``partial $G$-module'', where $G$ is a group,  is being understood in two related but different senses: initially the concept of a partial $G$-space was considered in \cite{DZh}, which means a $\kpar \,G$-module, or equivalently,  a $\kappa $-module $M$ equipped with a partial representation $G \to {\rm End}_{\kappa } (M),$ and which was generalized to Hopf case in  \cite{ABV} to deal with partial Hopf representations; the other one   means a partial action of $G$ on a commutative semigroup (whose domains are assumed to be ideals), and is used to deal with 
partial cohomology in   \cite{BaMoTe}, \cite{DKh},  \cite{DKh2},  \cite{DKh3},  \cite{DKhS},  \cite{DoPaPi1},  \cite{DoPaPiRo},  
  \cite{DoSa} and  \cite{KennedySchafhauser}. Of course, it is also natural  to take the first meaning and define the corresponding cohomology of $G,$ as it has been recently done in  \cite{AAR}. More precisely,  denote by $B$ the commutative subalgebra of $  \kpar \,G$ generated by all idempotents $e_g = [g][g\m],$ $g\in G$ (this subalgebra was denoted above by $\A$, but  we prefer now to reserve the symbol  $\A$ for a general algebra). The
partial representation 
$G\ni g\mapsto [g] \in  \kpar \,G$  results by (\ref{pRep->pAc}) in a unital partial action  $\tau $ of $G$ on $B,$ 
which in its turn gives rise to a partial representation $\lambda: G \to {\rm End}_{\kappa } (B),$  by the standard formula
\begin{equation}\label{standardFormula}
\lambda (g): B\ni b \mapsto \tau_g(be_{g\m}) = [g] b [g\m] \in B.
\end{equation} The latter endows  $B$ with a structure of a   left $\kpar \,G$-module,  and one may consider    $\kpar \,G$ as  an  augmented ring (see \cite[Chapter VIII]{CartanEilenberg}) by means of the  left  $\kpar \,G$-module epimorphism $\varepsilon :    \kpar \,G \to B,$ defined by  $$s=[g_1][g_2] \ldots [g_n] \mapsto  e_{g_1} e_{g_1 g_2} \ldots e_{g_1 g_2 \ldots g_n}=s s\m.$$ Then one takes the cohomology of the augmented ring $\kpar \,G \to B $ with values in an arbitrary       $\kpar \,G$-module $M$   \cite{AAR}, i.e.  one defines 
$$H_{par}^n(G,M)= {\rm  Ext}^n_{\kpar \,G} (B, M).$$  Notice the difference between this use of the notation $H_{par}^n$ and that in the above mentioned paper  \cite{BaMoTe}. Observe also that the ordinary group co\-ho\-mo\-lo\-gy can be seen as the cohomology of the augmented ring $\kappa G\to \kappa, g \mapsto 1,$ with values in a $\kappa G$-module, whereas,
the Hochschild cohomology of an associative algebra $\A $ can be defined as the cohomology of the augmented ring 
$\A \otimes \A ^{o} \to \A, a\otimes a' \mapsto a a',$ with values in an     $\A \otimes \A ^{o}$-module (i.e. an $\A$-bimodule), where $ \A ^{o}$ stands for the opposite algebra of $\A $ (see \cite{CartanEilenberg}). Given a unital partial action $\alpha $ of a group $G$ on an algebra $\A $ over a field $\kappa ,$ the authors  of   \cite{AAR} use Grothendieck's result \cite[Theorem 10.47]{Rotman} in order to prove the existence of a spectral sequence relating the 
 Hochschild cohomology of  the partial skew group ring  $\A\rtimes _\alpha G$ with   $H_{par}^n(G, - )$ and the
 Hochschild cohomology of $\A .$\\ 

To compare the cohomology from  \cite{AAR} with that from  \cite{DKh},  one may adapt  the projective resolution given in   \cite{DKh} to the  case of vector spaces, so that the $n$-cocycles and $n$-coboundaries are defined    by similar formulas. This may give the impression that the two cohomologies coincide.  Nevertheless, one should notice a significant difference in the concepts of $G$-modules. A   $\kappa G$-module
can be seen as a partial action of $G$ on a vector space   \cite[Example 5.4]{Ba}. In the definition of a partial group action  on a vector space (see  \cite[Remark 2.3]{AbDoExSi}) one assumes that each domain is a subspace and the involved bijections are $\kappa $-maps.    It is evident that the category of $\kappa G$-modules is abelian.  On the other hand, the partial $G$-modules from   \cite{DKh} are partial actions of $G$ on commutative semigroups. Their domains are unital two-sided ideals in the semigroup. Thus, if the semigroup under the partial action turns out to be a group, then the action has to be global. It is pointed out in    \cite{DKh} that the category of partial $G$-modules is not abelian. The requirement that the domains are ideals is essential to fit the theory of twisted partial group actions on algebras. As a consequence,  it allows one to define partial crossed products in the category of inverse semigroups  and interpret low-dimensional cohomology groups in terms of extensions. Furthermore, important structural results in semigroup theory use partial actions on meat-semilattices whose domains are order ideals. Semilattices can be considered as semigroups, and order ideals are semigroup ideals.\\

A natural situation in which both cohomologies appear is that of a unital partial action $\alpha $ of a group $G$ on a commutative algebra $\A $.  The domains are ideals generated by idempotents $1_g,$ $g\in G,$ in particular, they  are  vector subspaces. In this case  the  above  formula \eqref{standardFormula} becomes 
$$\pi (g): \A\ni a \mapsto \alpha_g(a 1_{g\m}) \in \A,$$ which gives  a partial representation $\pi $  of $G$ by 
$\kappa$-linear transformations  of $\A.$ Thus the cohomology from \cite{AAR} deals with the additive structure of $\A.$ Focusing on the multiplicative semigroup of $\A$ one has a partial $G$-module in the sense of   \cite{DKh} and the corresponding cohomology is based on the multiplicative structure of $\A.$ This shows that  the  cohomology groups from the two theories  have little in common, except for the case $n=0,$ in which both of them are related to partial in\-va\-ri\-ants.\\

The importance of groupoids for partial representations was  mentioned above. On the other hand, in \cite{AbadieThree} a   locally compact groupoid (called transformation groupoid) was related to any continuous partial action $\theta $ of a second countable locally compact Hausdorff group $G$ on a second countable locally compact Hausdorff  space $X$ and used to prove that the $C^*$-algebra of the groupoid is isomorphic to the  $C^*$-crossed product $C_0(X) \rtimes _{\theta}^{\rm full} G$ (here the  partial action of $G$ on $C_0(X),$ dual to that of $G$ on $X,$ is denoted by the same symbol $\theta $).
In the discrete group case another way to associate a groupoid to a partial action was given in \cite{AbadieFour}. More recent  results involving groupoids related to partial actions   include the groupoid $C^*$-algebra characterization of the enveloping  $C^*$-algebras  associated to   partial actions of  countable discrete groups on  locally compact spaces with applications in \cite{EGG}, as well as the study 
of the amenability of   groupoids related to  partial semigroup actions, englobing groupoids associated to   directed graphs, higher rank graphs and  topological higher rank graphs  in \cite{RenWil}. Moreover,  an algebraic version of the above mentioned  result from   \cite{AbadieThree} was established in \cite{BeuGon2} which says that given a continuous partial action of a discrete group $G$ on a locally compact totally disconnected  Hausdorff  space with clopen domains, the Steinberg algebra of the associated  transformation groupoid and the skew group ring by the corresponding  partial action of $G$ on the compactly supported locally constant function algebra are isomorphic. In the same article \cite{BeuGon2} the authors show that any Steinberg algebra (\cite{ClarkFarthingSimsTomforde}, \cite{St2010}), associated to an ample Hausdorff  groupoid, can be seen as a partial skew inverse semigroup ring.\\

The use of groupoids was essential in  \cite{MilanSt} for partial crossed product characterizations of  certain inverse semigroup $C^*$-algebras. This was done by dealing with  the above mentioned transformation groupoid of a partial group action from   \cite{AbadieThree} (also defined independently in  \cite{KL} in the discrete case), as well as   with the  universal groupoid ${\mathcal G} (S)$ of an inverse semigroup $S.$ The latter is the groupoid of germs of a certain natural action  of $S$ on the space of the semi-characters $\hat{E}$, where $E=E(S)$ stands for the idempotents of $S.$ The groupoid of germs of an inverse semigroup action is defined in a similar spirit with  the transformation groupoid of a partial action,  taking into account the natural partial order of $S.$ The reader is referred to 
 \cite{E4} for a detailed construction. The authors of \cite{MilanSt}  prove that if $S$ is a 
countable  $E$-unitary (or strongly $0$-$E$-unitary) inverse semigroups, then there is a partial action $\theta $ of the maximum group image $G$ of $S$ on  $\hat{E}$ such that  ${\mathcal G} (S)$ is isomorphic to the transformation groupoid of $\theta .$ This led  to the above mentioned crossed product descriptions. Moreover, the tight spectrum  $\hat{E}_{tight}$ (see  \cite{E4}) is invariant under the action of $S,$ and it is shown in  \cite{MilanSt}
that the corresponding groupoid of germs  ${\mathcal G}_{tight} (S)$ is also isomorphic to a transformation groupoid of a partial action of $G.$ As a consequence, the full $C^*$-algebra of      ${\mathcal G}_{tight} (S)$ is isomorphic to the crossed product  $C_0(  \hat{E}_{tight}) \rtimes ^{\rm full} G.$ The latter furnishes a general explanation for a number  of  partial  crossed product facts.\\

Apart from their importance to partial (semi)group actions and partial group representations, groupoids became an object of study also subject to their own partial actions. This was influenced by  the paper  \cite{Gil}, in which  the definition of a partial action of an ordered groupoid on a set was given, and it was proved that any such partial action is globalizable. Furthermore the Birget-Rhodes expansion $\tilde{G}^R$ of an ordered groupoid $G$ was defined and it was shown that there is a one-to-one correspondence between the partial actions of $G$ on a set $X$ and the actions of  $\tilde{G}^R$ on $X,$ extending  R.~Exel's result from \cite{E1}, mentioned above.\\ 

The study of partial groupoid action on rings began in the working group leaded by Antonio~Paques, the first publication being   \cite{BFP}, in which    the notion of a partial action $\alpha $ of an ordered groupoid $G$ on a ring $\A$ was introduced and  the corresponding partial skew groupoid ring  
$\A  \rtimes  _{\alpha}G $  was defined.  It was also proved that the meet-preserving actions of $\tilde{G}^R$ on $\A$ are
 in  a one-to-one correspondence with those  partial actions of  $G$ on  $\A$, in which the domain of each partial bijection is an ideal in $\A$. Furthermore, if $\alpha $ is a partial action of $G$
on $\A$ and $\beta $ is the corresponding action of $\tilde{G}^R$ on $\A$ then $\A  \rtimes  _{\alpha}G $ is an epimorphic image of  $\A  \rtimes  _{\beta} \tilde{G}^R, $ which extends a result from \cite{EV}, obtained in the case of  partial group actions. Analogous facts for partial actions of inductive groupoids on algebras were obtained in  \cite{B}, 
motivated by the fact that the category of inverse semigroup is  isomorphic to that of the inductive groupoids.\\

 The concept 
was further generalized in \cite{BP}, where the  definition of a partial action of an arbitrary groupoid on a ring was given, a criterion for its globalization was proved under the assumption that the identities of the groupoid act (trivially) on unital rings, showing also that the partial skew groupoid ring and the corresponding global one 
are Morita equivalent.   The latter fact was used to establish Galois theoretic results. More general globalization and Morita equivalence facts were  
recently obtained in \cite{BPi}, where the unital condition on the above mentioned domains was replaced by the $s$-unital condition. This extends also the results in \cite{DdRS} for partial group actions on $s$-unital rings.\\

The general definition from \cite{BP} gives, in particular, a notion of a groupoid action,  used  in   \cite{FlP}, where    Cohen-Montgomery-type duality theorems  for (co)actions of groupoids were proved, and in  \cite{FlP2} to study the Azumaya property for skew groupoid rings, obtaining  a nice  relation with  the notions of a Galois extension and Hirata-separability. 
  Another recent development is  a Galois theory in the style  of Grothendieck and Dress for groupoids acting on commutative rings \cite{PaTa}, including a Galois correspondence result.  The definition of a (global) groupoid action  in 
\cite{BP}   is   equivalent, under some conditions, to that given in   \cite{CaenDGr2}. Furthermore, in a latest article \cite{BPi2} necessary and sufficient conditions for the skew groupoid ring $A \ast _{\alpha } G$  to be a separable extension of $A$ were given, where $G$ is   a finite connected  groupoid acting partially on  a ring $A.$ For a general  $G$  with finite number of objects the problem of separability is reduced to the connected case. If an arbitrary groupoid $G$ acts globally on $A$ then the separability of the extension forces $G$ to be finite.\\

The globalization problem in the general setting of partial category actions  was investigated in the new article \cite{Nystedt}, in which the concept of a partial action of a category on a set is defined,  as well as that of a continuous partial action of a topological category on a topological space. It is proved that any partial action of a category on an arbitrary set admits a universal globalization, and the latter fact is also true in the context of continuous partial  category actions.  It is also shown that the new concept fits nicely the notion of a partial action of a groupoid on a set.\\

A new advance produced by the group of  A.~Paques comes from the fact that groupoid algebras  and Hopf algebras are both   englobed by   weak Hopf algebras.  The partial actions of weak Hopf algebras on rings were defined in
\cite{CasPaqQuaSant} and studied with respect to the globalization problem, and  the existence of a  Morita context  between the partial smash product and the corresponding global one.   A  difficulty  arises here when defining partial smash 
products  due to the fact that the tensor product is taken not over the ground field.  Amongst various  interesting  results a complete characterization of the partial actions of a weak Hopf algebra on the base field was given, suggesting a way to produce new examples. The authors also study  the relation between partial groupoid actions, in the sense of \cite{BP}, and partial actions of groupoid algebras as weak Hopf algebras,  showing how partial groupoid actions  fit into this general context.\\

A Hopf algebra $H$ may act partially on an algebra $\A$ from the right and from the left, and if these two actions commute, then we have the   concept of a partial $H$-bimodule algebra as defined by A.~Paques et al. in  \cite{CasPaqQuaSant2}. Similarly, the authors of  \cite{CasPaqQuaSant2} define partial   $H$-bicomodule algebras and study the globalization problem for both structures. 
The relations  between the partial $H^0$-biactions and  the partial $H$-bicoactions is studied, where $H^0$ stands for the finite dual of $H.$ In addition,  under the assumption that  $H^0$ separates points, the authors prove that  the globalization of a partial $H^0$-bimodule algebra and  that of the corresponding  partial $H$-bicomodule algebra are isomorphic as $H^0$-bimodule algebras.\\

As a dual counterpart of the globalization theorems for partial Hopf (co)actions on algebras in \cite{AB2} and \cite{AB3},  the globalization problem for partial  Hopf  (co)actions on coalgebras was recently considered in \cite{CasQua}. More precisely,
it is proved that every partial module coalgebra is globalizable and, under   a certain rationality condition, the globalization for a partial comodule coalgebra  also exists. The relations between  actions and coactions on algebras and coalgebras are analyzed, as well as between their globalizations (sometimes assuming some natural conditions).  Interesting examples are also given.\\

Multiplier algebras entered the theory of partial actions at the stage of its foundations: R.~Exel's definition of a continuous twisted 
partial action of a locally compact group on a $C^*$-algebra uses  invertible multipliers of certain ideals   as values of the twist \cite{E0}. This was transferred to  abstract algebras  \cite{DES1}, as well as to semigroups   \cite{DKh3}. In purely ring theoretic setting, multipliers turned out to be  crucial to deal with the as\-so\-ci\-a\-ti\-vi\-ty of partial crossed pro\-ducts \cite{DE}, \cite{DES1}, whose arguments were  also used for semigroups in  \cite{DKh}, \cite{DKh2} and \cite{DKh3}, and which also provided  an algebraic proof of the as\-so\-ci\-a\-ti\-vi\-ty of partial $C^*$-crossed products, which was initially achieved by means of approximate identities.   Furthermore, multipliers are essential ingredients of the crossed product criterion  for graded algebras \cite[Theorem 6.1]{DES1}, as well as of the more general partial cohomology theory  \cite{DKh3}.\\

 On the other hand, multipliers came into the Hopf theoretic picture in A.~Van~Daele's definition of ``non-unital Hopf algebras'', i.e. the so-called multiplier Hopf algebras    \cite{VanDaele}, in which the co-multiplication  takes values in the multiplier algebra of the tensor pro\-duct, $\Delta : H \to M(H \otimes H).$ This concept was born in 1994, and since then much attention is being payed to various aspects of multiplier Hopf algebras and their generalizations. In particular, A.~Paques introduced the authors of 
the hot of print article \cite{AzCaFoMa} to the subject and suggested to investigate their actions and co-actions on non-necessarily unital algebras. The notion of a partial action of a regular multiplier Hopf algebra $H$ on an algebra $R$ with a nondegenerate product is introduced and studied in  \cite{AzCaFoMa}, as well as that of a coaction of $H$ on $R.$ The relation between the actions and coactions is discussed. The main result gives, under natural assumptions,  
 a Morita context connecting the coinvariant algebra $R^{\underline{coH}}$ with an appropriate   smash product. Finally 
 a notion of a partial Galois coaction is related  to this Morita context.\\

  A rather general version of Hopf Galois theory was elaborated in \cite{AlFeGoSo}, dea\-ling with the so-called lax  entwining structures  introduced by S.~Caenepeel and K.~Janssen in  \cite{CaenJan0}\footnote{The preprint \cite{CaenJan0} was a preliminary version of \cite{CaenJan}.}, and their particular cases of weak and partial  en\-t\-wi\-ning structures. Weak entwining structures can be obtained  from comodule algebras over weak bialgebras, whereas  partial entwining structures can be constructed from partial actions of bialgebras on algebras.  Each of these  entwining structures consists of an algebra $\A$ and a coalgebra ${\mathcal C}$ together with a map $\psi : {\mathcal C}\otimes \A \to \A \otimes {\mathcal C}$ obeying certain axioms. A classical result by Y.~Doi and M.~Takeuchi  (1986) states that an $H$-extension  is cleft precisely  when it is Hopf Galois with normal basis.  The authors extend this fact to the general context of lax Galois extensions in a symmetric monoidal category, by introducing  lax ${\mathcal C}$-cleft extensions and proving that the lax ${\mathcal C}$-cleft extensions are exactly the lax ${\mathcal C}$-Galois extensions with normal basis. Further results on  ${\mathcal C}$-Galois extensions, as well as on various types of   entwining structures were obtained in  
\cite{Soneira2016}.\\

An important part of the research on partial actions is dedicated to  ideals and the  ring theoretic properties of partial crossed products.
First of all, one asks whether or not the partial crossed products are associative.  It turned out that it is not always the case. In \cite{DE} a ring $\A$ 
was called {\it strongly associative} if for any group $G$ and any partial action $\alpha $ of $G$ on $\A,$ the partial skew group ring $\A\rtimes _\alpha G$ is associative. In the   language of multipliers a sufficient condition for the associativity  of 
$\A \rtimes _\alpha G$ was given in \cite{DE},  in particular, semiprime rings are strongly associative. This was extended in \cite{BFP} to partial skew groupoid rings. On the other hand, the 
ring of the upper-triangular $n\times n$-matrices over a field is strongly associative if and only if $n\leq 2.$ Another example of the non-associativity  was given by constructing a partial action of the group of order $2$ on a four-dimensional algebra, which in the case of characteristic $2$ is isomorphic to the modular group algebra of the Klein-four group. This suggested the problem of the characterization of strongly associative modular group algebras of finite groups, which was solved in 
\cite{Jonas}, \cite{Jonas2} and \cite{Jonas3}. In the recent preprint \cite{Jonas5} the strongly associative modular semigroup algebras $\kappa S$ are characterized, where   $S$ is a finite inverse semigroup  without zero, which is not a group, and $\kappa $ is an algebraically closed field whose characteristic $p$ divides the order of the kernel of $S.$ As a consequence, it is pointed out that the modular partial group algebra $\kpar \, G$ is never strongly  associative, provided that $G$ is a finite group and $\kappa $ is an algebraically closed field whose characteristic divides $|G|.$\\

 A systematic  ring-theoretic  investigation  of  partial skew group rings,  the rings under  partial actions and the rings under their enveloping actions   began in \cite{FL}, whereas the study of the prime and maximal ideals in \cite{CF} initiated the consideration  of  partial skew polynomial rings.    Subsequent results were obtained in   \cite{CFM}, \cite{AvLa}, \cite{AvFeLa},  \cite{C}, \cite{CFHM}, \cite{Avi1}, \cite{AvFe}, \cite{CCF}, \cite{C2}  and \cite{CFG}. The treatment of the ring-theoretic properties of the more general structure of a crossed product by a twisted partial group action began in  \cite{BLP} and \cite{PSantA}. In particular,   results on the Azumaya property from  \cite{PSantA} were important to define in \cite{DoPaPi1}   a homomorphism from the partial cohomology group 
$H^2(G,\af, \A )$ to the relative Brauer group  $B(\A /{\A }^\af)$ by sending any arbitrary partial $2$-cocycle  to the corresponding crossed product.  This map fits into  the seven terms exact sequence \cite{DoPaPiRo}, mentioned above, and in the  case of a (global) Galois extension of fields results in the classical isomorphism between the second cohomology group and the relative Brauer group.\\ 

More recent  developments on the ring theoretic properties include the study of the simplicity of the partial crossed products
in  \cite{Gon}, \cite{GonOinRoy},  \cite{BaraCoSo}, \cite{NysOin-1}  and \cite{NysOin}, and the chacarterization of the Leavitt path algebras as partial group rings \cite{GR}, which permitted one to obtain  alternative proofs for the simplicity criterion and for the   
Cuntz-Krieger uniqueness theorem for Leavitt path algebras (see \cite{GR} and \cite{GonOinRoy}).  Furthermore,  in  \cite{AvLa3} the relations between the module properties over $\A ,$ $\A \rtimes _{\alpha} G$ and $\A ^{\alpha}$ were studied, where   $\A ^{\alpha}$ denotes the subring of the $\alpha$-invariants, whereas  in \cite{CoGo} the related structure of the partial skew 
power series ring was recently considered with respect to the distributive and  Bezout duo properties.\\ 

As  to  more general skew structures,   in   \cite{GonYon} the Leavitt path algebras were  characterized as partial skew groupoid rings, which was applied to  study a class of groupoid graded isomorphisms between  Leavitt path algebras.
In \cite{BeuGonOinRoy}  the simplicity of the skew group ring by  a partial action of an   inverse semigroup on  a commutative ring was characterized  and used to offer a new proof for the simplicity  criterion 
 for a Steinberg algebra associated with a Hausdorff and ample groupoid. Furthermore, in 
   \cite{CavalStAna1}  the se\-mi\-pri\-mi\-ti\-vi\-ty and the se\-mi\-pri\-ma\-li\-ty properties  for partial smash products were considered, as well as their prime and Jacobson radicals.\\

In a recent preprint  \cite{NysOinPin}  a rather general situation of  a unital partial (in particular,  global) action of a groupoid $G$ on a non-necessarily associative ring $R$ was considered. In particular, it was proved that if the partial action is unital and $R$ is alternative, then  $R\ast G$ is left (or right) artinian if and only if so too is $R$ and $R_g=\{0\}$ for all but 
 finitely many   $g\in {\rm mor } ( G).$ In the global case the assumption on $R$ to be alternative can be dropped. This extends the well-known result  of 1963 by I.~G.~Connel
on group rings, as well as  its more general version for global skew group rings  by J.~K.~Park  (1979). Applications were given to  groupoid rings, partial group rings, generalized matrix rings and Leavitt path algebras, and, moreover,  the noetherian pro\-per\-ty for partial skew groupoid rings was also considered.\\

 For the twisted case  in \cite{BeCoFeFlo} the authors use the above discussed weak globalization for a twisted partial action $\alpha $ of a group $G$ on a semiprime ring to give  several results   with necessary and sufficient conditions for a partial crossed product $\A  \ast _{\alpha} G$ to be a semiprime left Goldie ring, imposing rea\-so\-nable restrictions, such as $G$ being finite and $\A $ having no $|G|$-torsion, or $G$ being an infinite polycyclic group and $\alpha $ being of finite type etc.  
In order to obtain these results the authors establish  a series  of facts relating   the  pro\-per\-ties  and ideals of 
$\A \ast _{\alpha} G$ to those of $Q_m(\A ) \ast _{\alpha} G.$ In particular, if  $\A  \ast _{\alpha} G$ is semiprime left Goldie ring, then so too is $\A ,$ without any restriction on $G$ and $\alpha.$\\ 
 
More recently,  in  \cite{BeCo} partial crossed product by twisted partial actions were considered with respect to the artinian, noetherian, semilocal, perfect and   semiprimary    properties, under some assumptions on the twisted partial actions, such as to be of finite type. If the twisted partial action is globalizable, then the above ring theoretic properties for the  partial crossed products were related to those for the corresponding global crossed product.
 Moreover, the Krull dimension  of partial crossed products was also studied  in \cite{BeCo}, as well as the global dimension and the weak global dimension. Furthermore, triangular matrix representations of partial skew group rings were also considered, and,  in addition, it was shown with certain assumptions that the ring under the twisted partial action is Frobenius or symmetric, then so too is the crossed product.\\  

In another preprint \cite{CoRu}  twisted partial skew power series rings and twisted partial skew Laurent series rings are introduced and  studied with respect to the properties of being    a prime ring,  a  semi-prime ring, a semiprime Goldie ring or a Goldie ring with semiprime base subring. Prime ideals are considered as well as the prime radical (the latter for the twisted partial skew power series case).\\

In relation to the twisted case in \cite{NysOinPin2} the authors introduce the concept of a unital ring  $S$  which is epsilon strongly graded by a group $G,$ as a si\-mul\-ta\-ne\-o\-us generalization of the strongly graded rings and crossed products by unital twisted 
partial actions.  They define a trace map and using it give a criterion of the separability of the extension $S_1 \subseteq S.$ This generalizes a known result for strongly graded rings, as well as  the corresponding  result from the above mentioned paper \cite{BLP}. Amongst other facts,  the hereditary and Frobenius properties, and the simplicity  for epsilon-strongly graded rings are also   discussed, and   the so-called  epsilon-crossed products are  introduced, which are proved to  be isomorphic  with the  unital partial crossed products.  A series of examples are worked out, including those  giving  separable epsilon-strongly graded rings, which are neither  strongly graded, nor partial crossed products in some natural way.\\

As to the Hopf case, in \cite{CavalStAna2} a  general theory of partial $H$-radicals is de\-ve\-loped   for partial $H$-module algebras based on the Amitsur-Kurosh general radical theory. Using the monoidal category of left partial $H$-modules from \cite{ABV} and its subcategory of algebra objects, a definition of a partial smash product 
$\underline{A \#  H}$  is elaborated, where $H$ is a Hopf algebra with a bijective antipode and $A$ is a non-necessarily unital $H$-module algebra (in the above mentioned subcategory).    Several examples of  partial $H$-radicals of Jacobson type are worked out and stu\-di\-ed. The obtained results extend  facts established  in the global case in \cite{Fisher}   and \cite{Sidorov}.  Some questions posed in \cite{Fisher} are also discussed.\\

In \cite{Zhang2016} for a finite dimensional Hopf algebra $H$ over a field, imposing a condition of a  lazy $1$-cocycle  \cite{BichonKassel},
the concept of a strong partial $H$-module algebra is defined.  
For strong partial $H$-module algebras the partial trace map is considered and a Maschke-type theorem is proved.  Se\-pa\-ra\-bi\-li\-ty of the partial smash product over the base ring is also discussed.\\ 

The partial trace map naturally occurred in earlier papers, and it is related to both Galois theory and Morita theory. Amongst the above mentioned  Morita theoretic facts, two type of Morita equivalence results are  used in the algebraic  investigations of ring theoretic properties. The first one has its origin in \cite{AbadieTwo}, where it was proved that if $\alpha $ is a globalizable partial action of a locally compact Hausdorff group $G$ on a $C^*$-algebra $\A $ and  $(\B, \beta ) $ is its enveloping action, then the reduced $C^*$-crossed products $\A \rtimes _{\alpha}^{\rm red}   G$ and    $\A \rtimes _{\beta}^{\rm red}   G$ are strongly Morita equivalent. We already  referred to the recent $C^*$-theoretic extensions of this fact in \cite{AbBussFerraro}  and \cite{AbFerraro2}.  Algebraic versions   appeared in \cite{DE}, \cite{DdRS}, \cite{DES1},  \cite{AB2}, \cite{BP}, \cite{AAB}, \cite{CasPaqQuaSant}, \cite{BPi}. The other type of results is inspired by the global case and for  partial group actions   appeared  in   \cite{AvLa}:  given a unital partial action of a group $G$ on a ring $\A ,$ there is a natural Morita context linking  the subring of the $\alpha $-fixed points $\A ^{\alpha}$ and the skew group ring $\A \rtimes _{\alpha}   G;$  the context is strict if and only if $\A ^{\alpha} \subseteq \A $ is an $\alpha $-Galois extension. Its analogue for partial Hopf actions was established in \cite{AB1}, and  Hopf theoretic  considerations about  the corresponding Morita context  is the topic of a recent  survey by A.~Paques \cite{Paques2016}. The discussion in  \cite{Paques2016}  involves the partial trace map, the strictness of the context, Galois extensions, separability,   and is accompanied    by a series of   illustrative examples.\\

Considerations on ideals in partial crossed products naturally apper when dealing with  ring-theoretic properties, however, 
they are of great interest  on their own. Thus in one of the latest papers \cite{C3} the author investigate the prime Goldie 
ideals in the related structure of the partial skew polynomial rings, whereas in   \cite{BeuGon} the ideals are studied in 
skew group rings by partial actions on algebras of functions on an abstract set $X$.\\

 More precisely, let $\kappa $ be a field
and denote by ${\mathcal F}_0(X)$ the algebra of all functions $X\to \kappa $ with finite support, endowed with the pointwise operations.  Then there is a bijection between the non-zero ideals of  ${\mathcal F}_0(X)$ and the non-empty subsets of $X,$ and, moreover, there is a one-to-one correspondence between the partial actions of a group $G$ on $X$ and the partial actions of $G$ on ${\mathcal F}_0(X).$   If the partial action $\alpha $ of $G$ on the set $X$ is free, then the skew group ring  
${\mathcal F}_0(X) \rtimes _{\alpha} G$ is shown to be isomorphic to a certain equivalence relation algebra, and this is 
used to describe  the ideals of   ${\mathcal F}_0(X) \rtimes _{\alpha} G$ as function rings on appropriate subsets of $X.$
The results are  analogues of  known facts on  partial actions on the  $C^*$-algebra $C_0(X),$ where $X$ is a locally compact Hausdorff space. As a motivation, the authors include a proof of a simplified version of  the above mentioned 
  F. Abadie's characterization of $C^*$-partial crossed products as groupoid $C^*$-algebras, for the case of free partial actions of countable groups acting on unital commutative $C^*$-algebras.\\

There is a rich history of the study of ideals in $C^*$-algebraic crossed pro\-ducts, with a special role played by the so-called
 Effros-Hahn conjecture  \cite{EffrosHahn}, which in its original form states that every primitive ideal in the crossed product of a commutative $C^*$-algebra by a locally compact group should be induced from a primitive ideal in the $C^*$-algebra of some isotropy group. It was proved in \cite{Sauvageot} in the  case of discrete amenable groups and later extended to various  contexts.\\

 In  \cite{DE3} the  Effros-Hahn conjecture is studied in the algebraic setting, more specifically, for the skew group rings    ${\mathscr L}_c (X) \rtimes G,$ where $G$ is an arbitrary discrete group acting partially on a locally compact Hausdorff,
totally disconnected  topological space $X,$ and ${\mathscr L}_c (X)$ stands for  the algebra  of all locally constant, compactly supported functions  $X$ with values in an arbitrary  field $\kappa .$ It is also assumed that each domain $X_g$ in the partial action is  closed and open (clopen)  in $X.$ One of the main results in \cite{DE3} says that every   ideal in ${\mathscr L}_c (X) \rtimes G$ is an intersection of  ideals induced from isotropy groups.\\

 Tools are developed in order to  understand the induction process, one of them being the notion of an {\it admissible} ideal of the  group algebra $\kappa H_{x_0},$  where $H_{x_0}$ is the isotropy group of a point $x_0 \in X.$ In particular,  denoting  by 
${\rm Ind}_{x_0}(I)$ the ideal of    ${\mathscr L}_c (X) \rtimes G$ induced from  an ideal of $I$ of $\kappa H_{x_0},$ the map $I \mapsto  {\rm Ind}_{x_0}(I)$ is  proved to be an injection of the set of the  admissible ideals of  $\kappa H_{x_0}$ into the set of the ideals of     ${\mathscr L}_c (X) \rtimes G.$ In order to compare ideals induced from distinct isotropy groups the  notion of  a {\it transposition} of ideals
is introduced, which takes  any ideal of $ \kappa H_{x_0}$ to an admissible ideal  of $\kappa H_{\hat{x}_0},$ ${x}_0, \hat{x}_0 \in X.$ It is proved in \cite{DE3} that   given admissible ideals $I \lhd \kappa H_{x_0}$ and $\hat{I} \lhd \kappa H_{\hat{x}_0},$ the equality  
${\rm Ind}_{x_0}(I) =   {\rm Ind}_{\hat{x}_0}(\hat{I})$ occurs  if and only if $\hat{I}$ is the transposition of $I.$ A concrete description of a transposed ideal is given, and  it is also shown  that the induction preserves some important properties of the ideals, namely, being primitive, or prime, or 
meet-irreducible.\\

Notice that the above mentioned assumption on $X_g$ to be closed (besides being open by definition)  naturally appears 
in \cite{EGG}: the enveloping space of a topological partial action of a countable discrete group on a compact Hausdorff space is Hausdorff  if and only if each domain $X_g$ is closed (note that we already mentioned that this condition was also used in \cite{BeuGon2}).\\

The ideas of  \cite{DE3} were  further developed in the most recent preprint \cite{Demeneghi} to the study of the ideals in the   Steinberg algebra 
$A_{\kappa }({\mathcal G})$ of an ample Hausdorff groupoid $G$ over a field $\kappa$, by considering  $A_{\kappa }({\mathcal G})$ as as an inverse semigroup crossed product algebra. As in \cite{DE3}, it is proved that   every ideal in   $A_{\kappa }({\mathcal G})$ may be obtained as the intersection of ideals induced from isotropy group algebras. Since the partial skew group ring    ${\mathscr L}_c (X) \rtimes G$ from   \cite{DE3} can be viewed as 
the Steinberg algebra \cite{St2010}  for the transformation groupoid associated to the partial action of $G$ on $X$, the  algebraic version of the Effros-Hahn conjecture proved in \cite{DE3} is extended now to the generality of   Steinberg algebras.\\

Recent works  around partial actions include also
 the study 
of the category of the partial Doi-Hopf modules in \cite{ChenWang2013},
 of partial actions on  power sets     in \cite{AvLa2}, 
 of the category of partial $G$-sets for a fixed group $G$ in    \cite{AviBuiZap},
of partial orbits and $n$-transitivity  in   \cite{AvHe-GoOr-Ja}, 
of  partial group entwining structures and   partial group (co)actions of a Hopf group coalgebra on a family of algebras in \cite{ChenWang2014},
 of generalized partial smash products in \cite{DuLin}, and
  of twisted partial Hopf coactions and corresponding   partial crossed coproducts in 
\cite{ChenWangKang2017}, as well as  a note on sums of ideals  \cite{AvLa2011}. More information around partial actions may be found  
 in the surveys \cite{Ba}, \cite{D}, \cite{D2}, \cite{F2}, \cite{Paques2016}, \cite{Paques2018} and \cite{Pi4}.  Notice also that partial  actions and/or related structures have been  mentioned in \cite{AlFeGo}, \cite{AlFeGo2}, \cite{AraE2},  \cite{BatiCaenVerc}, \cite{bohm}, \cite{BuE2011}, \cite{ChLi0}, \cite{ChLi2006},  \cite{DoomsVeloso}, \cite{E4}, \cite{ExelPardo}, \cite{FerGonRod2016} and \cite{FerGonRod-ArXiv2015}.\\

 \section*{Acknowledgements} 

The author  thanks Fernando Abadie and Mykola Khrypchenko  for many useful comments.


\bibliographystyle{amsplain}


\end{document}